\documentclass[a4paper,12pt]{amsart}
\usepackage[utf8]{inputenc}
\usepackage[english]{babel}
\usepackage[text={160mm,230mm},centering]{geometry}
\usepackage{amsmath,amssymb,amsthm,mathabx,verbatim}
\usepackage{amssymb,latexsym}
\usepackage{tikz-cd}
\usepackage[bookmarks=true]{hyperref}
\usepackage{booktabs}
\usepackage{caption}
\usepackage{enumerate}
\usepackage{cleveref}

\numberwithin{equation}{section}
\numberwithin{table}{section}

\newcommand{\tx}{\textnormal}
\newcommand{\QQ}{\mathbb Q}
\newcommand{\ZZ}{\mathbb{Z}}
\newcommand{\CC}{\mathbb{C}}
\newcommand{\RR}{\mathbb R}
\newcommand{\FF}{\mathbb F}
\newcommand{\HH}{\mathbb{H}}
\newcommand{\PP}{\mathbb P}

\newcommand{\End}{\mathrm{End}}

\newcommand{\Norm}{\mathrm{Norm}}
\newcommand{\Aut}{\mathrm{Aut}}
\newcommand{\GL}{\mathrm{GL}}
\newcommand{\Frob}{\mathrm{Frob}}
\newcommand{\GLn}{\mathrm{GL}_2(\mathbb Z/n\mathbb Z)}

\newcommand{\SLZ}{\mathrm{SL}_2(\mathbb Z)}

\newcommand{\new}{\textnormal{new}}
\newcommand{\ns}{\textnormal{ns}}
\newcommand{\s}{\textnormal{s}}

\newcommand{\smt}[4]{\left( \begin{smallmatrix} #1 &#2\\ #3 &#4\\ \end{smallmatrix} \right) }

\newtheorem{thm}[equation]{Theorem}%[section] %teorema numerato
\newtheorem{prop}[equation]{Proposition} %proposizione numerata
\newtheorem{lem}[equation]{Lemma} %lemma numerato
\newtheorem{Algorithm}[equation]{Algorithm}
\newtheoremstyle{named}{}{}{\itshape}{}{\bfseries}{.}{.5em}{#3}
\theoremstyle{named} 
\theoremstyle{remark}
\newtheorem{rem}[equation]{Remark} %remark numerato
\newtheorem{ese}[equation]{Example} %esempio numerato
\theoremstyle{definition}
\title{Modular Curves with many points over Finite Fields}
\author{Valerio Dose}
\address{Valerio Dose, Department of Computer, Control and Management Engineering, ``Sapienza'' University of Rome, Roma, Italy}
\email{valerio.dose@uniroma1.it}
\author{Guido Lido}
\address{Guido Lido, Department of Mathematics, University of Rome ``Tor Vergata'', Roma, Italy}
\email{guidomaria.lido@gmail.com}
\thanks{The second author is supported by the MIUR ``Excellence Department Project MATH@TOV'',
awarded to the Department of Mathematics, University of Rome ``Tor
Vergata'', and by the ``Programma Operativo (PON) “Ricerca e Innovazione”
2014-2020'' } 
\author{Pietro Mercuri}
\address{Pietro Mercuri, Dipartimento SBAI, ``Sapienza'' University of Rome, Roma, Italy}
\email{mercuri.ptr@gmail.com}
\author{Claudio Stirpe}
\address{Claudio Stirpe, Convitto Nazionale ``R. Margherita'', Anagni, Italy}
\email{clast@inwind.it}


\begin{document}
%\oddsidemargin=-2.0 cm
%\evensidemargin=-2.0 cm
%\textwidth
%\topmargin=-2.5 cm
%\textheight 28.2 cm
\maketitle

\begin{abstract}
%We give an algorithm to compute the number of points over finite fields of modular curves $X_H$, associated to subgroups $H$ of $\GL_2(\mathbb Z/n\mathbb Z)$, where for each prime $p$ dividing $n$, the subgroup $H$ at $p$ is either a Borel subroup, a Cartan subgroup, or the normalizer of a Cartan subgroup of  $\GL_2(\mathbb Z/p^e\mathbb Z)$. We also consider quotients of such curves by Atkin-Lehner involutions $w_p$ whenever $H$ is a Borel subgroup at $p$. 
%We apply our algorithm to curves with small genus and to fields with small cardinality, improving many lower bounds for the maximum number of points over a finite field of a curve with fixed genus. To make the computation we generalize Chen's isogeny to modular curves of composite level and mixed Borel, split/non-split Cartan type.

We describe an algorithm to compute the number of points over finite fields on a broad class of modular curves: we consider quotients $X_H/W$ for $H$ a subgroup of $\GL_2(\mathbb Z/n\mathbb Z)$ such that for each prime $p$ dividing $n$, the subgroup $H$ at $p$ is either a Borel subroup, a Cartan subgroup, or the normalizer of a Cartan subgroup of  $\GL_2(\mathbb Z/p^e\mathbb Z)$, and for $W$ any subgroup of the Atkin-Lehner involutions of $X_H$. We applied our algorithm to more than ten thousands curves of genus up to 50, finding more than one hundred record-breaking curves, namely curves $X/\FF_q$ with genus $g$ that improve the previously known lower bound for the maximum number of points over $\FF_q$ of a curve with genus $g$. As a key technical tool for our computations, we prove the generalization of Chen's isogeny to all the Cartan modular curves of composite level.
\end{abstract}

\section{Introduction}

Finding among the algebraic curves of a  fixed genus, the ones with the largest number of points over a finite field of fixed cardinality, is an interesting effort in algebraic geometry and number theory which also has applications in coding theory (see for example \cite[Section 8.4]{Sti09}). The Weil bound prescribes the inequality
\[
\#C(\mathbb{F}_q)\leq q+1+ 2g\sqrt q,
\]
where $\#C(\mathbb{F}_q)$ is the number of points over a finite field $\FF_q$, with $q$ being a prime power, of a nonsingular, projective, absolutely irreducible curve $C$ of genus $g$. 

Let $q$ be a fixed prime power. As it is shown in \cite{DV83} and \cite{XS95}, the previous estimate cannot be sharp for $g$ large since
$$
N_g(\mathbb F_q)\le g(\sqrt{q}-1)+o(g),\quad \text{for }g\to +\infty,
$$
where $\displaystyle N_g(\mathbb F_q):=\max_{C\text{ of genus }g}\#C(\mathbb F_q)$. Several sequences $\{C_n\}_{n\in\mathbb N}$ of algebraic curves over $\mathbb F_q$ with increasing genus $g_n$, have been found to achieve the asymptotic bound $\displaystyle\lim_{n\to+\infty}\frac{\#C_n(\mathbb F_q)}{g_n}=\sqrt{q}-1$. Among these sequences, one of the most classical
example is obtained by taking modular curves $X_0(m)$ with certain increasing levels $m$ and counting supersingular points over a field $\mathbb F_q$, where $q$ is a square, see \cite{TVZ82}. Concerning curves with small genus, the website \cite{manypoints} is devoted to collect the full list of known upper bounds $M_g(\mathbb{F}_q)$ and lower bounds $L_g(\mathbb{F}_q)$ for $N_g(\mathbb F_q)$, when $g\leq 50$ and for fields of characteristic less than $100$  with small cardinality. 
%We denote by $M_g(\mathbb{F}_q)$ and $L_g(\mathbb{F}_q)$, respectively, these upper and lower bounds for the maximum number of points of a genus $g$ curve over $\mathbb{F}_q$. 
Notice that $M_g(\mathbb{F}_q)$ and $L_g(\mathbb{F}_q)$ depend on the current state of the art in this line of research.

In this paper we look at this question concentrating on a large class of modular curves. Namely, we consider modular curves $X_H$ associated to a subgroup $H$ of $\GL_2(\mathbb Z/n\mathbb Z)$ such that, for each prime $p$ dividing $n$, with maximum power $e$, the reduction of $H$ modulo $p^e$ is either a Borel subgroup, a Cartan subgroup, or the normalizer of a Cartan subgroup of  $\GL_2(\mathbb Z/p^e\mathbb Z)$. We also consider quotients of these curves by Atkin-Lehner involutions $w_p$ whenever $H$ is a Borel subgroup at $p$. We give an algorithm which computes the number of points over a finite field on these modular curves, without having the equation of the curve, but only using the trace of Hecke operators acting on modular abelian varieties associated to weight 2 newforms invariant under the congruence subgroup $\Gamma_0(m)$. Numerical approximations of these traces are available at \cite{lmfdb} for $m\le 10000$. We applied our algorithm to these modular curves over small finite fields of characteristic less than 20 using data available on \cite{lmfdb}. Furthermore, in September 2022 we compared our results with the best known curves available on the database \cite{manypoints} and we found many improvements that we list in the final section of this paper.
%old: We computed the number of points over small finite fields of characteristic less than 20 for every curve $X_H/\langle w_p\colon p|n_{0^+}\rangle$ with $H$ a subgroup of $\GL_2(\mathbb Z/n\mathbb Z)$ for $n=n_0n_{0^+}n_\ns n_{\ns^+}$ where $n_0,n_{0^+},n_\ns,n_{\ns^+}$ are pairwise coprime and $H$ is Borel at each prime dividing $n_0n_{0^+}$, non-split Cartan at each prime divinding $n_\ns$, normalizer of a non-split Cartan at each prime dividing $n_{\ns^+}$, and for $n_0n_{0^+}n_\ns^2 n_{\ns^+}^2\le 10000$.
% C: We computed the number of points over finite fields of characteristic less than 20 for every curve $X_H/\langle w_p\colon p|n_{0^+}\rangle$ 
% with $H$ a subgroup of $\GL_2(\mathbb Z/n\mathbb Z)$ 
% with $n=n_0n_{0^+}n_\ns n_{\ns^+}$ 
% where $n_0,n_{0^+},n_\ns,n_{\ns^+}$ are pairwise coprime with $n_0n_{0^+}n_\ns^2 n_{\ns^+}^2\le 10000$.
% We use the notation:
% \begin{itemize}
% \item $H$ is 
% Borel at each prime dividing $n_0n_{0^+}$;
% \item $H$ is non-split Cartan at each prime divinding $n_\ns$; 
% \item $H$ is normalizer of a non-split Cartan at each prime dividing $n_{\ns^+}$.
% \end{itemize}

To make the computation, we prove that the Jacobian of  the modular curves $X_H$ we are considering, are isogenous to a product of modular abelian varieties associated to weight 2 newforms invariant under the congruence subgroup $\Gamma_0(m)$. We give explicitly this factorization of the Jacobian (\Cref{th:mixedjac}), generalizing results of \cite{Che98}, \cite{dSE00}, \cite{Che04}, \cite{DLM22}.

The idea to use these types of modular curves builds on recent work on computing equations, rational points, and automorphism groups for such curves, as for example \cite{Mer18}, \cite{MS20}, \cite{DMS19}, \cite{DLM22}, \cite{Dos16}, \cite{DFGS14}, \cite{M-R22}, \cite{BDMTV19}, \cite{BFG21}, \cite{frengley2023congruences}, \cite{AABCCKW21}. Another algorithm has been recently devised for computing the number of points on $X_H$ for a general $H$, see \cite[Section 5]{RSZ22}. Their method does not make use of the factorization of the Jacobian of the modular curve up to isogeny, but they rather can obtain it case by case as a consequence of the computation of the number of points of the curves over many fields of different characteristic \cite[Section 6]{RSZ22}.

The paper is organized as follows. In \Cref{sec:not} we introduce the definitions related to the modular curves we are considering. 
%In \Cref{sec:supersingular} we recall the theory about the number of supersingular points on modular curves and we discuss the related asymptotic estimates .
%In \Cref{sec:jac} we prove the theorem which shows that the Jacobian of our modular curves is isogenous to an explicit product of modular abelian varieties, each of them associated to some weight 2 newform invariant under the congruence subgroup $\Gamma_0(d)$ for some level $d$. 
In \Cref{sec:jac} we describe the Jacobians of our curves, up to isogeny, explicitly in terms of the Jacobian of Borel modular curves. In \Cref{sec:basic} we discuss the algorithm we use to compute the number of points over finite fields. 

In \Cref{sec:asympt} we give asymptotic results to estimate the number of $\FF_q$-points on our curves when the genus $g$ is large. Finally, in \Cref{sec:app} and in the Appendix, we collect the results obtained. In particular, for every choice of $g$ and $\FF_q$, we list in the Appendix the curve with the largest number of points among the ones we considered, and in \Cref{sec:app}, for the convenience of the reader, we collect all the examples which improve the previously known lower bounds $L_g(\FF_q)$.

\section{Modular curves of mixed type}\label{sec:not}

%Let $n$ be a positive integer and $\Gamma$ be a subgroup of $\SL_2(\ZZ)$. Such groups act on the complex upper half-plane $\HH=\{ \tau \in \CC: \mathrm{Im}(\tau)>0 \}$, see \cite{DS05}. We denote by $\Gamma(n)$ the subgroup of $\SL_2(\ZZ)$ given by the matrices $\left(\begin{array}{c c} a&b\\c&d \end{array}\right),$ with entries $b\equiv c\equiv 0 \bmod n$ and $a\equiv d\equiv 1 \bmod n$. In this section we consider a special family $\Gamma_0(n)$ of subgroups of $\SL_2(\ZZ)$ containing $\Gamma(n)$, namely matrices $\left(\begin{array}{c c} a&b\\c&d \end{array}\right),$ with $c\equiv 0 \bmod n$. We denote by $X_0(n)$ the modular curve associated to $\Gamma_0(n)\backslash \HH^*$, where $\HH^*=\HH\cup\PP^1(\QQ)$.

Let $n$ be a positive integer. For each subgroup $H$ of $\GLn$, we define
\[
\Gamma_{H}:= \{ \gamma \in \tx{SL}_2(\ZZ): \gamma^T\!\!\!\pmod n \tx{ lies in }H\}.
\]
Let $\HH=\{ \tau \in \CC: \mathrm{Im}(\tau)>0 \}$ be the complex upper half-plane and denote by $\HH^*=\HH\cup\PP^1(\QQ)$ the extended complex upper half-plane. We define the modular curve associated to $H$:
\[
X_H:=\Gamma_H\backslash \HH^*.
\]
For a more detailed reference about this construction see \cite[Section 1]{DLM22}. In the following, we choose a non-square element $\xi\in (\ZZ/p^e\ZZ)^\times$ when $p$ is an odd prime and $e$ is a positive integer. 
We define the following subgroups of $\textnormal{GL}_2(\ZZ/p^e\ZZ)$ for every prime $p$:
\allowdisplaybreaks
\begin{align*}
&B^0(p^e):=\left\{\begin{pmatrix}a & 0 \\ c & d \end{pmatrix}, a,c,d\in\ZZ/p^e\ZZ, \quad ad \not\equiv 0 \bmod p \right\}; \\
&C_{\textnormal{s}}(p^e):=\left\{\begin{pmatrix}a & 0 \\ 0 & d \end{pmatrix}, a,d\in(\ZZ/p^e\ZZ)^\times \right\}; \\
&C_{\textnormal{s}}^+(p^e):=C_{\textnormal{s}}(p^e)\cup\left\{\begin{pmatrix}0 & b \\ c & 0 \end{pmatrix}, b,c\in(\ZZ/p^e\ZZ)^\times \right\}; \\
&C_\ns(2^e):=	\left\{\begin{pmatrix}a & b \\ b & a+b \end{pmatrix}, a,b\in\ZZ/2^e\ZZ,  (a,b) \not\equiv (0,0) \bmod 2 \right\}; \\
&C_\ns^+(2^e):=C_\ns(2^e)\cup\left\{\begin{pmatrix}a & a-b \\ b & -a \end{pmatrix}, a,b\in\ZZ/2^e\ZZ, (a,b) \not\equiv (0,0) \bmod 2 \right\}; \\
&C_\ns(p^e):=	\left\{\begin{pmatrix}a & b\xi \\ b & a \end{pmatrix}, a,b\in\ZZ/p^e\ZZ,  (a,b) \not\equiv (0,0) \bmod p \right\}, \quad \text{if $p$ is odd}; \\ 
&C_\ns^+(p^e):=C_\ns(p^e)\cup\left\{\begin{pmatrix}a & b\xi \\ -b & -a \end{pmatrix}, a,b\in\ZZ/p^e\ZZ, (a,b) \not\equiv (0,0) \bmod p \right\}, \quad \text{if $p$ is odd}; \\
&B_r(p^e):=\left\{\begin{pmatrix}a & bp^r \\ cp^{r+1} & d \end{pmatrix}, a,b,c,d\in\ZZ/p^e\ZZ, \quad ad \not\equiv 0 \bmod p \right\}, \quad \text{for }r=0,1,\ldots,e-1; \\
&T_r(p^e):=\left\{\begin{pmatrix}a & bp^r\\ cp^r & d \end{pmatrix}, a,b,c,d\in\ZZ/p^e\ZZ, ad-bcp^{2r} \in (\ZZ/p^e\ZZ)^\times \right\}.
\end{align*}
We remark that $T_e(p^e)=C_{\textnormal{s}}(p^e)$ and that $C_{\textnormal{s}}(p^e), C_\ns(p^e)$ are respectively a split and a non-split Cartan subgroup of $\GL_2(\ZZ/p^e\ZZ)$ and, with the exception of $C^+_{\textnormal{s}}(2^e)$, the groups $C_{\textnormal{s}}^+(p^e), C_\ns^+(p^e)$ are the corresponding normalizers inside $\GL_2(\ZZ/p^e\ZZ)$.
\allowdisplaybreaks[0]

Given  $n_0$, $n_{0^+}$, $n_{\textnormal{s}}$, $n_{\textnormal{s}^+}$, $n_{\textnormal{ns}}$, $n_{\textnormal{ns}^+}$  pairwise coprime positive integers such that $n=n_0n_{0^+}n_{\textnormal{s}}n_{\textnormal{s}^+}n_{\textnormal{ns}}n_{\textnormal{ns}^+}$. By Chinese Remainder Theorem we have $\GLn\cong \prod_{i=1}^r \GL_2(\ZZ/p_i^{e_i}\ZZ)$. We look at the subgroups of $\GLn$ of the following form
\begin{equation} \label{eq:H=prod_H_j}
H \cong \prod_{i=1}^r H_{p_i}, \quad \textnormal{where } H_{p_i}=\begin{cases}
B^0(p_i^{e_i}), & \text{if }p_i\mid n_0n_{0^+}, \\
C_{\textnormal{s}}(p_i^{e_i}), & \text{if }p_i\mid n_{\textnormal{s}}, \\
C_{\textnormal{s}}^+(p_i^{e_i}), & \text{if }p_i\mid n_{\textnormal{s}}^+, \\
C_\ns(p_i^{e_i}), & \text{if }p_i\mid n_{\textnormal{ns}}, \\
C_\ns^+(p_i^{e_i}), & \text{if }p_i\mid n_{\textnormal{ns}}^+. \\
\end{cases}
\end{equation}
Then we define our modular curves of mixed type as
\begin{align}
&X(n_0,n_{0^+}, n_{\textnormal{s}}, n_{\textnormal{s}^+}, n_{\textnormal{ns}}, n_{\textnormal{ns}^+}):=X_H/\langle w_{p_i},\text{ for every prime $p_i$ dividing }n_{0^+}\rangle,\label{eq:modular-curve-symbol} \\
&X(n_0,n_{0^+}, n_{\textnormal{ns}}, n_{\textnormal{ns}^+}):=X(n_0,n_{0^+},1,1, n_{\textnormal{ns}}, n_{\textnormal{ns}^+}),
\label{eq:modular-curve-1-symbol}
\end{align}
where $w_{p_i}$ denotes the Atkin-Lehner operator associated to $p_i$. Let $n=p_1^{e_1}\cdots p_r^{e_r}$ be the prime factorization of $n$.
We also define, for every positive integer $n$, the congruence subgroup
\[
\Gamma_0(n):=\{\left(\begin{smallmatrix}a&b\\ c&d\end{smallmatrix}\right)\in\mathrm{SL}_2(\ZZ):\left(\begin{smallmatrix}a&b\\ c&d\end{smallmatrix}\right)\equiv\left(\begin{smallmatrix}*&*\\ 0&*\end{smallmatrix}\right) \mod n\},
\]
we denote the $\CC$-vector space of the cusp forms of weight 2 invariant under $\Gamma_0(n)$ by $\mathcal S_2(\Gamma_0(n))$ and by $\mathcal S_2^\new(\Gamma_0(n))$ the subspace generated by the newforms. We define the modular curve:
\[
X_0(n):=\Gamma_0(n)\backslash \HH^*.
\]
We also denote the Jacobian of this curve by $J_0(n):=\mathrm{Jac}(X_0(n))$ and its new part by $J_0^\new(n)$, which is the factor of $J_0(n)$ isogenous to the abelian variety associated to $\mathcal S_2^\new(\Gamma_0(n))$.

\begin{rem}\label{rem:s_non_serve}

There is an isomorphism of algebraic curves 
$$
X(n_0, n_{0^+}, n_\s, n_{\s^+}, n_\ns, n_{\ns^+}) \cong X(n_0 n_\s^2, n_{0^+} n_{\s^+}^2, n_\ns, n_{\ns^+}).
$$
Indeed, writing the above curves respectively as $X_{H_1}/W_1 = \Gamma_1\backslash \HH^*$ and $X_{H_2}/W_2 = \Gamma_2\backslash \HH^*$, for $\Gamma_1, \Gamma_2$ subgroups of $\mathrm{PGL}_{2}^{\det>0}(\RR) = \Aut(\HH)$, the group $\Gamma_1$ can be obtained conjugating $\Gamma_2$ firstly by $\smt{0}{-1}{n_\s n_{\s^+}}{0}$ and then conjugating again by a suitable matrix in $\mathrm{SL}_2(\ZZ)$. 

%and either $a\equiv d\not\equiv 0 \pmod{n_0n_{0^+}n_\ns n_{\ns^+}}$ if $(n_\s n_{\s^+})^{-3}$ is a square modulo $n_0n_{0^+}n_\ns n_{\ns^+}$ or $a\equiv \xi d\not\equiv 0 \pmod{n_0n_{0^+}n_\ns n_{\ns^+}}$, for $\xi$ a non-square modulo $n_0n_{0^+}n_\ns n_{\ns^+}$, if $(n_\s n_{\s^+})^{-3}$ is not a square modulo $n_0n_{0^+}n_\ns n_{\ns^+}$
%G: bisogna anche imporre condizioni su $c,d$? P: seguono da quelle che ho scritto mi pare quando mi ero posto il problema, c dovrebbe avere la stessa condizione di b e d la stessa di a; comunque se preferisci mi va anche bene la prima versione che hai scritto di seguito G: Altrimenti, in questo nostro cruccio infinito  propongo anche $M$ tale che $M \equiv \smt{n_\s^2 n_{\s^+}^2(1+n_0n_{0^+}n_\ns n_{\ns^+})}{-n_0n_{0^+}n_\ns n_{\ns^+}}{n_\s n_{\s^+}n_0n_{0^+}n_\ns n_{\ns^+}}{n_\s^2 n_{\s^+}^2(1+n_0n_{0^+}n_\ns n_{\ns^+})}$ modulo $ n_\s^2 n_{\s^+}^2n_0n_{0^+}n_\ns n_{\ns^+}$ e con determinante uguale a $n_\s n_{\s^+}$. Oppure $M$ congruo a $ \smt{0}{-1}{n_\s n_{\s^+}}{0}$ modulo $ n_\s^2 n_{\s^+}^2$, congruo all'identita' modulo $n_0n_{0^+}n_\ns n_{\ns^+}$ e con determinante $n_\s n_{\s^+}$. V: Bellissimo, ma forse troppo per un lettore che lo legge la prima volta, ovvero solleva più domande di quante ne risponde. Forse lascerei una versione glissante, e poi in caso rompiamo il culo (e i coglioni) al referee in fase di review.)

\end{rem}

\section{Jacobians of modular curves of mixed type and generalization of Chen's isogeny} \label{sec:jac}

With the notation introduced in Section \ref{sec:not}, in this section we show that the Jacobian of the curves $X(n_0,n_{0^+}, n_{\textnormal{s}}, n_{\textnormal{s}^+}, n_{\textnormal{ns}}, n_{\textnormal{ns}^+})$ is isogenous to a product of Atkin-Lehner quotients of $J_0^\new(m)$ for suitable levels $m$. We start with the key lemma where we discuss separately the four Cartan cases.
\begin{lem}\label{lem:special_cases}
Let $n>1$ be an integer and let $H< \GL_2(\ZZ/n\ZZ)$ be a subgroup. We use the following notation:
\[
J_0^\new(n)^{w_{p_1}w_{p_2}\ldots w_{p_k}}:=J_0^\new(n)/\langle w_{p_1},w_{p_2},\ldots, w_{p_k}\rangle,
\]
for $p_1,\ldots,p_k$ distinct primes dividing $n$ and $w_{p_j}$ the Atkin-Lehner involution associated to $p_j$, for $j=1,\ldots,k$. Then we have:
\begin{enumerate}
\item\label{lem:special_cases_s} If $H=C_{\textnormal{s}}(n):= \prod_{i=1}^r C_{\textnormal{s}}(p_i^{e_i})$, then
\begin{equation}
\mathrm{Jac}(X_H)=J_{\textnormal{s}}(n) \sim \prod_{d|n^2}J_0^\new(d)^{\sigma_0\left(\frac{n^2}{d}\right)},
\end{equation}
where $\sigma_0(m)$ is the number of positive divisors of an integer $m$.
\item\label{lem:special_cases_s+} If $H=C_{\textnormal{s}}^+(n)  := \prod_{i=1}^r C_{\textnormal{s}}^+(p_i^{e_i})$, then
\begin{equation}
\mathrm{Jac}(X_H)=J_{\textnormal{s}}^+(n)  \sim \prod_{d|n^2}J_0^\new(d)^{\sigma_0^*\left(\frac{n^2}{d}\right)},
\end{equation}
where $\sigma_0^*$ is the function defined by
\[
\sigma_0^*(p^f):=\begin{cases}
\frac{1}{2}(f+1), & \text{if }f\text{ is odd}, \\
\frac{f}{2}+w_p, & \text{if }f\text{ is even},
\end{cases}
\]
for a prime $p$ and a positive integer $f$, and $\sigma_0^*(m_1m_2)=\sigma_0^*(m_1)\sigma_0^*(m_2)$, for $m_1,m_2$ coprime positive integers.
\item\label{lem:special_cases_ns} If $H=C_{\textnormal{ns}}(n) := \prod_{i=1}^r C_{\textnormal{ns}}(p_i^{e_i})$, then
\begin{equation}
\mathrm{Jac}(X_H)=J_{\textnormal{ns}}(n)\sim \prod_{d|n}J_0^\new(d^2).
\end{equation}
\item\label{lem:special_cases_ns+} If $H=C_{\textnormal{ns}}^+(n)  := \prod_{i=1}^r C_{\textnormal{ns}}^+(p_i^{e_i})$, then
\begin{equation}
\mathrm{Jac}(X_H)=J_{\textnormal{ns}}^+(n) \sim \prod_{d|n}J_0^\new(d^2)^{w_{p_1}\ldots w_{p_s}},
\end{equation}
where $p_1,\ldots, p_s$ are all the primes dividing $d$.
\end{enumerate}
\end{lem}

\begin{proof}
Part \ref{lem:special_cases_s}. We have
\begin{equation*}
\mathrm{Jac}(X_H)=J_{\textnormal{s}}(n) \cong J_0(n^2),
\end{equation*}
because $W_n^{-1}\Gamma_{C_\textnormal{s}(n)}(n)W_n=\Gamma_0(n^2)$, where $W_n:=\smt{0}{-1}{n}{0}$, and
\begin{equation*}
J_0(n^2) \sim\prod_{d|n^2}J_0^\new(d)^{\sigma_0\left(\frac{n^2}{d}\right)},
\end{equation*}
because of classical Atkin-Lehner theory about the oldforms (see \cite{DS05}).

Part \ref{lem:special_cases_s+}. Let $n=p_1^{e_1}\ldots p_k^{e_k}$ be the prime factorization of $n$ and let $W_Q$, for a positive integer $Q$, be the matrix defined in \cite[after Lemma~7]{AL70}. We have
\begin{equation*}
\mathrm{Jac}(X_H)=J_{\textnormal{s}}^+(n)  \cong J_0(n^2)^{w_{p_1}\ldots w_{p_k}},
\end{equation*}
because $\smt{0}{-1}{n}{0}\Gamma_{C_\textnormal{s}^+(n)}\smt{0}{-1}{n}{0}^{-1}$ is congruent to $\left\langle\Gamma_0(n^2),W_{p_1^{2e_1}},\ldots ,W_{p_k^{2e_k}}\right\rangle$ modulo scalar matrices.
Then
\begin{equation*}
J_0(n^2)^{w_{p_1}\ldots w_{p_k}}\sim\prod_{d|n^2}J_0^\new(d)^{\sigma_0^*\left(\frac{n^2}{d}\right)},
\end{equation*}
because of Atkin-Lehner theory (see \cite[Equations 5.1 and 5.2]{AL70}).

Part \ref{lem:special_cases_ns}.  
Let $n=p_1^{e_1}\ldots p_k^{e_k}$ be the prime factorization. For each $c=p_1^{f_1}\ldots p_k^{f_k}$, we define 
\[
K(c):=\prod_{j=1}^k K_j(p_j^{f_j}) < \GL_2(\ZZ/c\ZZ), \quad \text{with} \quad K_j(p_j^{f_j}):=\begin{cases}
T_{\frac{f_j}{2}}(p_j^{e_j}), & \text{if }f_j \text{ is even},\\
B_{\frac{f_j-1}{2}}(p_j^{e_j}), & \text{if } f_j\text{ is odd}.
\end{cases}
\]
Then, using the machinery in \cite{dSE00} together with \cite[Proposition 3.2]{DLM22}, we get
\begin{equation*}
\mathrm{Jac}(X_H)=J_{\textnormal{ns}}(n) \sim \prod_{c|n^2}\mathrm{Jac}(X_{K(c)})^{\varepsilon(c)m(c)},
\end{equation*}
%because of Chen's isogeny described in \cite[Section~3]{DLM22}, 
where the functions $\varepsilon(c)$ and $m(c)$ are defined by
\begin{equation}\label{eq:eps-m}
\varepsilon(p^f):=(-1)^f, \quad\text{and}\quad m(p^f):=\begin{cases}
1, & \text{if }p^f||n^2, \\
2, & \text{otherwise},
\end{cases}
\end{equation}
for a prime power dividing $n^2$ and by $\varepsilon(d_1d_2)=\varepsilon(d_1)\varepsilon(d_2)$ and $m(d_1d_2)=m(d_1)m(d_2)$ for $d_1,d_2$ coprime divisors of $n^2$. (For example if $n^2=p^2q^2$ and $c=p=p^1q^0$, then $m(c)=m(p^1)m(q^0)=4$.) Moreover, we have
\begin{equation*}
 \prod_{c|n^2}\mathrm{Jac}(X_{K(c)})^{\varepsilon(c)m(c)} \cong \prod_{c|n^2}J_0(c)^{\varepsilon(c)m(c)},
\end{equation*}
because $\Gamma_{K(c)}$ and $\Gamma_0(c)$ are conjugate, as explained in \cite[proof of Theorem~3.8]{DLM22}. Then we have 
\begin{equation*}
 \prod_{c|n^2}J_0(c)^{\varepsilon(c)m(c)}\sim\prod_{c|n^2}\prod_{d|c}J_0^\new(d)^{\varepsilon(c)m(c)\sigma_0\left(\frac{c}{d}\right)},
\end{equation*}
because of classical Atkin-Lehner theory about the oldforms (see \cite{DS05}). Last equality:
\begin{equation*}
\prod_{c|n^2}\prod_{d|c}J_0^\new(d)^{\varepsilon(c)m(c)\sigma_0\left(\frac{c}{d}\right)}= \prod_{d|n}J_0^\new(d^2),
\end{equation*}
follows by
\begin{align*}
\prod_{c|n^2}\prod_{d|c}J_0^\new(d)^{\varepsilon(c)m(c)\sigma_0\left(\frac{c}{d}\right)}&=\prod_{d|n^2}\prod_{d|c|n^2}J_0^\new(d)^{\varepsilon(c)m(c)\sigma_0\left(\frac{c}{d}\right)}= \\
&=\prod_{d|n^2}J_0^\new(d)^{\sum_{d|c|n^2}\varepsilon(c)m(c)\sigma_0\left(\frac{c}{d}\right)};
\end{align*}
and, if $n=p_1^{e_1}\ldots p_k^{e_k}$ and $d=p_1^{g_1}\ldots p_k^{g_k}$ and $c=p_1^{f_1}\ldots p_k^{f_k}$ are the prime factorizations of $n,d,c$, we have
\begin{align*}
&\sum_{d|c|n^2}\varepsilon(c)m(c)\sigma_0\left(\frac{c}{d}\right)=\sum_{f_1=g_1}^{2e_1}\ldots\sum_{f_k=g_k}^{2e_k}\varepsilon(p_1^{f_1}\ldots p_k^{f_k})m(p_1^{f_1}\ldots p_k^{f_k})\sigma_0\left(p_1^{f_1-g_1}\ldots p_k^{f_k-g_k}\right)=\\
&=\prod_{j=1}^k\sum_{f_j=g_j}^{2e_j}\varepsilon(p_j^{f_j})m(p_j^{f_j})\sigma_0\left(p_j^{f_j-g_j}\right)=\prod_{j=1}^k\left(2e_j-g_j+1+\sum_{f_j=g_j}^{2e_j-1}(-1)^{f_j}2(f_j-g_j+1)\right)=\\
&=\prod_{j=1}^k\left(2e_j-g_j+1+2\sum_{f_j=g_j}^{2e_j-1}(-1)^{f_j}f_j+2(1-g_j)\sum_{f_j=g_j}^{2e_j-1}(-1)^{f_j}\right)=\prod_{j=1}^k\left(\frac{1}{2}(1+(-1)^{g_j})\right)=\\
&=\begin{cases}
1, & \text{if }g_1\equiv\ldots\equiv g_k\equiv 0 \bmod 2, \\
0, & \text{otherwise},
\end{cases}
\end{align*}
where in the fifth equality we used
\[
\sum_{i=a}^{b-1} x^i=\frac{x^a-x^b}{1-x}\quad \text{and}\quad\sum_{i=a}^{b-1} ix^{i-1}=\frac{x^a-x^b}{(1-x)^2}+\frac{ax^{a-1}-bx^{b-1}}{1-x},
\]
hence
\[
\prod_{d|n^2}J_0^\new(d)^{\sum_{d|c|n^2}\varepsilon(c)m(c)\sigma_0\left(\frac{c}{d}\right)}=\prod_{d|n}J_0^\new(d^2).
\]

Part \ref{lem:special_cases_ns+}. 
Given the factorization $n=p_1^{e_1}\ldots p_k^{e_k}$, for each $c=p_1^{f_1}\ldots p_k^{f_k}$, we define 
\[
K'(c)\!:=\!\prod_{j=1}^k K_j'(p_j^{f_j}) \!<\! \GL_2(\ZZ/c\ZZ), \,\, \text{with} \,\, K'_j(p_j^{f_j})\!:=\!\begin{cases}
T_{\frac{f_j}{2}}(p_j^{e_j}), & \!\!\text{if }f_j \ne 2e_j \text{ and even},\\
C_{\textnormal{s}}^+(p_j^{e_j}), & \!\!\text{if }f_j =2e_j,\\
B_{\frac{f_j-1}{2}}(p_j^{e_j}), & \!\!\text{if } f_j\text{ is odd}.
\end{cases}
\]
Then, using the machinery explained in \cite{dSE00}, together with \cite[Theorem 1.1]{Che04} and its extension to the even case described in Section~3 and in the appendix of \cite{DLM22}, we get 
\begin{equation*}
\mathrm{Jac}(X_H)=J_{\textnormal{ns}}^+(n) \sim\prod_{c|n^2}\mathrm{Jac}(X_{K'(c)}) ^{\varepsilon(c)},
\end{equation*}
where $\varepsilon$ is defined as in the proof of Part~\ref{lem:special_cases_ns}.

Then, defining the matrices $W_Q$, for a positive integer $Q$, as in \cite[after Lemma~7]{AL70} 
and since $\smt{0}{-1}{m}{0}\Gamma_{K'(c)}\smt{0}{-1}{m}{0}^{-1}$ is congruent to $\left\langle\Gamma_0(c),W_{p_j^{f_j}} \text{ for }p_j | c \text{ s.t. }f_j=2e_j\right\rangle$ modulo scalar matrices, where $m$ is such that $c=m^2\prod_{\substack{j=1 \\ f_j \text{ odd}}}^k p_j$, we have
\begin{equation*}
\prod_{c|n^2}\mathrm{Jac}(X_{K'(c)}) ^{\varepsilon(c)}\cong\prod_{c|n^2}J_0(c)^{\varepsilon(c)\prod'w_{p_j}},
\end{equation*}
where
\[
\textstyle\prod' w_{p_j}:=\displaystyle\prod_{\substack{p_j | c \text{ s.t.} \\ f_j=2e_j}}w_{p_j}.
\]
Then
\begin{equation*}
\prod_{c|n^2}J_0(c)^{\varepsilon(c)\prod'w_{p_j}}\sim\prod_{c|n^2}\prod_{d|c}J_0^\new(d)^{\varepsilon(c)\sigma_0'\left(\frac{c}{d}\right)},
\end{equation*}
is similar to Parts \ref{lem:special_cases_s} and \ref{lem:special_cases_s+}, where
\begin{equation}\label{eq:sigma'}
\sigma_0'(p^f):=\begin{cases}
\sigma_0^*(p^f), & \text{if } p^f||\frac{n^2}{d}, \\
\sigma_0(p^f), & \text{otherwise},
\end{cases}
\end{equation}
for a prime power dividing $\frac{n^2}{d}$ and $\sigma_0'(d_1d_2)=\sigma_0'(d_1)\sigma_0'(d_2)$ for $d_1,d_2$ coprime divisors of $\frac{n^2}{d}$. Last equality:
\begin{equation*}
\prod_{c|n^2}\prod_{d|c}J_0^\new(d)^{\varepsilon(c)\sigma_0'\left(\frac{c}{d}\right)}= \prod_{d|n}J_0^\new(d^2)^{w_{p_1}\ldots w_{p_k}},
\end{equation*}
follows by
\begin{align*}
\prod_{c|n^2}\prod_{d|c}J_0^\new(d)^{\varepsilon(c)\sigma_0'\left(\frac{c}{d}\right)}&=\prod_{d|n^2}\prod_{d|c|n^2}J_0^\new(d)^{\varepsilon(c)\sigma_0'\left(\frac{c}{d}\right)}= \\
&=\prod_{d|n^2}J_0^\new(d)^{\sum_{d|c|n^2}\varepsilon(c)\sigma_0'\left(\frac{c}{d}\right)};
\end{align*}
and, if $d=p_1^{g_1}\ldots p_k^{g_k}$ is the prime factorization of $d$, we have
\begin{align*}
&\sum_{d|c|n^2}\varepsilon(c)\sigma_0'\left(\frac{c}{d}\right)=\sum_{f_1=g_1}^{2e_1}\ldots\sum_{f_k=g_k}^{2e_k}\varepsilon(p_1^{f_1}\ldots p_k^{f_k})\sigma_0'\left(p_1^{f_1-g_1}\ldots p_k^{f_k-g_k}\right)=\\
&=\prod_{j=1}^k\sum_{f_j=g_j}^{2e_j}\varepsilon(p_j^{f_j})\sigma_0'\left(p_j^{f_j-g_j}\right)=\prod_{j=1}^k\left(\sigma_0^*\left(p_j^{2e_j-g_j}\right)+\sum_{f_j=g_j}^{2e_j-1}(-1)^{f_j}(f_j-g_j+1)\right)=\\
&=\prod_{j=1}^k\left(\sigma_0^*\left(p_j^{2e_j-g_j}\right)+\sum_{f_j=g_j}^{2e_j-1}(-1)^{f_j}f_j+(1-g_j)\sum_{f_j=g_j}^{2e_j-1}(-1)^{f_j}\right)=\\
&=\prod_{j=1}^k\left(\sigma_0^*\left(p_j^{2e_j-g_j}\right)-e+\frac{1}{4}\Big((-1)^{g_j}-1+2g_j\Big)\right)=\\
&=\begin{cases}
\prod_{j=1}^k w_{p_j}, & \text{if }g_1\equiv\ldots\equiv g_k\equiv 0 \bmod 2, \\
0, & \text{otherwise},
\end{cases}
\end{align*}
where in the fifth equality we used
\[
\sum_{i=a}^{b-1} x^i=\frac{x^a-x^b}{1-x}\quad \text{and}\quad\sum_{i=a}^{b-1} ix^{i-1}=\frac{x^a-x^b}{(1-x)^2}+\frac{ax^{a-1}-bx^{b-1}}{1-x},
\]
hence
\[
\prod_{d|n^2}J_0^\new(d)^{\sum_{d|c|n^2}\varepsilon(c)\sigma_0'\left(\frac{c}{d}\right)}=\prod_{d|n}J_0^\new(d^2)^{w_{p_1}\ldots w_{p_k}}.
\]
\end{proof}
\begin{thm}\label{th:mixedjac}
Let $X:=X(n_0$, $n_{0^+}, n_{\textnormal{s}},n_{\textnormal{s}^+}, n_{\textnormal{ns}}, n_{\textnormal{ns}^+})$ be as in \Cref{eq:modular-curve-symbol}, then
\begin{align}\label{Cartan-Borel-isogeny-fine}
\mathrm{Jac}(X)\sim\prod_{d\mid N}J_0^\new(d_0d_{0^+}d_{\textnormal{s}} d_{\textnormal{s}^+} d_{\textnormal{ns}}^2 d_{\textnormal{ns}^+}^2)^{\tiny\sigma_0\bigg(\frac{n_0n_{\textnormal{s}}^2}{d_0d_{\textnormal{s}}}\bigg)\sigma_0^*\bigg(\frac{n_{0^+}n_{\textnormal{s}^+}^2}{d_{0^+}d_{\textnormal{s}^+}}\bigg) \displaystyle\prod_{p\in \mathcal P(d_{\textnormal{ns}^+})}w_p},
\end{align}
where $N:=n_0n_{0^+}n_{\textnormal{s}}^2n_{\textnormal{s}^+}^2n_{\textnormal{ns}}n_{\textnormal{ns}^+}$ and $d:=d_0d_{0^+}d_{\textnormal{s}} d_{\textnormal{s}^+} d_{\textnormal{ns}} d_{\textnormal{ns}^+}$, with $d_0\mid n_0$, $d_{0^+}|n_{0^+} $, $d_{\textnormal{s}}|n_{\textnormal{s}}^2$,  $d_{\textnormal{s}^+}|n_{\textnormal{s}^+}^2$, $d_{\textnormal{ns}}|n_{\textnormal{ns}}$, $d_{\textnormal{ns}^+}|n_{\textnormal{ns}^+}$, and
\[
J_0^\new(p_1^{e_1}p_2^{e_2}\ldots p_k^{e_k})^{w_{p_1}w_{p_2}\ldots w_{p_k}}:=J_0^\new(p_1^{e_1}p_2^{e_2}\ldots p_k^{e_k})/\langle w_{p_1},w_{p_2},\ldots, w_{p_k}\rangle,
\]
for $p_1,\ldots,p_k$ distinct primes and $w_{p_j}$ is the Atkin-Lehner involution associated to $p_j$, for $j=1,\ldots,k$; moreover $\mathcal P(m)$ is the set of prime divisors of an integer $m$ and $\sigma_0(m)$ is the number of positive divisors of $m$ and $\sigma_0^*$ is the function defined by
\[
\sigma_0^*(p^f):=\begin{cases}
\frac{1}{2}(f+1), & \text{if }f\text{ is odd}, \\
\frac{f}{2}+w_p, & \text{if }f\text{ is even},
\end{cases}
\]
for a prime $p$ and a positive integer $f$, and $\sigma_0^*(m_1m_2)=\sigma_0^*(m_1)\sigma_0^*(m_2)$, for $m_1,m_2$ coprime positive integers. 
\end{thm}

\begin{proof}
This is a natural extension of the previous lemma and the argument is the same. Given $c_{\textnormal{ns}}$ and $c_{\textnormal{ns}^+}$ divisors of $n_{\textnormal{ns}}^2$ and $n_{\textnormal{ns}^+}^2$ respectively, for each $c:=n_0n_{0^+}n_{\textnormal{s}}n_{\textnormal{s}^+}c_{\textnormal{ns}}c_{\textnormal{ns}^+}$, we define
\[
K(c):=\prod_{\substack{p^f||c\\p\text{ prime}}} K(p^f),
\]
where, if $e$ is the $p$-adic valuation of $n$, we define
\[
K(p^f):=\begin{cases}
B^0(p^{f}), & \text{if }p|n_0n_{0^+}, \\
C_{\textnormal{s}}(p^{f}), & \text{if }p|n_{\textnormal{s}}, \\
C_{\textnormal{s}}^+(p^{f}), & \text{if }p|n_{\textnormal{s}^+}, \\
B_{\frac{f-1}{2}}(p^{e}), & \text{if }p|c_{\textnormal{ns}}c_{\textnormal{ns}^+} \text{ and if } f\text{ is odd}, \\
T_{\frac{f}{2}}(p^{e}), & \text{if }p|c_{\textnormal{ns}}c_{\textnormal{ns}^+} \text{ and $f$ is even and }f\ne 2e, \\
C_{\textnormal{s}}(p^{e}), & \text{if }p|c_{\textnormal{ns}} \text{ and }f=2e, \\
C_{\textnormal{s}}^+(p^{e}), & \text{if }p|c_{\textnormal{ns}^+} \text{ and }f=2e.
\end{cases}
\]
Again, using the machinery explained in \cite{dSE00}, together with the representation theoretic results in \cite{Che04} and \cite{DLM22}, and defining $\varepsilon$, $m$, $\sigma_0'$ as in Equations (\ref{eq:eps-m}) and (\ref{eq:sigma'}), we find that
\begin{align*}
&\mathrm{Jac}(X_H)^{\prod w_{p_{0^+}}}\sim \\
&\sim \prod_{\substack{c_{\textnormal{ns}}|n_{\textnormal{ns}}^2\\c_{\textnormal{ns}^+}|n_{\textnormal{ns}^+}^2}} \mathrm{Jac}(X_{K(n_0n_{0^+}n_{\textnormal{s}}n_{\textnormal{s}^+}c_{\textnormal{ns}}c_{\textnormal{ns}^+})})^{\varepsilon(c_{\textnormal{ns}})\varepsilon(c_{\textnormal{ns}^+})m(c_{\textnormal{ns}})\prod w_{p_{0^+}}}\cong \\
&\cong \prod_{\substack{c_{\textnormal{ns}}|n_{\textnormal{ns}}^2\\c_{\textnormal{ns}^+}|n_{\textnormal{ns}^+}^2}} J_0(n_0n_{0^+}n_{\textnormal{s}}^2n_{\textnormal{s}^+}^2 c_{\textnormal{ns}}c_{\textnormal{ns}^+})^{\varepsilon(c_{\textnormal{ns}})\varepsilon(c_{\textnormal{ns}^+})m(c_{\textnormal{ns}})\prod w_{p_{0^+}}\prod w_{p_{\textnormal{s}^+}}\prod'w_{p_{\textnormal{ns}^+}}}\sim \\
&\sim \prod_{\substack{c_{\textnormal{ns}}|n_{\textnormal{ns}}^2\\c_{\textnormal{ns}^+}|n_{\textnormal{ns}^+}^2}} \prod_{d\mid N} J_0^{\mathrm{new}}(d)^{\varepsilon(c_{\textnormal{ns}})\varepsilon(c_{\textnormal{ns}^+})m(c_{\textnormal{ns}})\sigma_0\left(\frac{n_0}{d_0}\right)\sigma_0^*\left(\frac{n_{0^+}}{d_{0^+}}\right)\sigma_0\left(\frac{n_{\textnormal{s}}^2}{d_{\textnormal{s}}}\right)\sigma_0^*\left(\frac{n_{\textnormal{s}^+}^2}{d_{\textnormal{s}^+}}\right)\sigma_0\left(\frac{c_{\textnormal{ns}}}{d_{\textnormal{ns}}}\right)\sigma_0'\left(\frac{c_{\textnormal{ns}^+}}{d_{\textnormal{ns}^+}}\right)}= \\
&= \prod_{d\mid N} J_0^{\mathrm{new}}(d_0d_{0^+}d_{\textnormal{s}}d_{\textnormal{s}^+} d_{\textnormal{ns}}^2d_{\textnormal{ns}^+}^2)^{\sigma_0\left(\frac{n_0}{d_0}\right)\sigma_0^*\left(\frac{n_{0^+}}{d_{0^+}}\right)\sigma_0\left(\frac{n_{\textnormal{s}}^2}{d_{\textnormal{s}}}\right)\sigma_0^*\left(\frac{n_{\textnormal{s}^+}^2}{d_{\textnormal{s}^+}}\right)\prod w_{p_{\textnormal{ns}^+}}},
\end{align*}
where $d$ and $N$ are defined in the statement of the theorem and the product $\prod w_{p_{0^+}}$ varies over all the primes $p_{0^+}$ dividing $n_{0^+}$, the product $\prod w_{p_{\textnormal{s}^+}}$ varies over all the primes $p_{\textnormal{s}^+}$ dividing $c_{\textnormal{s}^+}$, the product $\prod w_{p_{\textnormal{ns}^+}}$ varies over all the primes $p_{\textnormal{ns}^+}$ dividing $n_{\textnormal{ns}^+}$ and the product $\prod'w_{p_{\textnormal{ns}^+}}$ varies over all the primes $p_{\textnormal{ns}^+}$ such that $v_{p_{\textnormal{ns}^+}}(c_{\textnormal{ns}^+})=v_{p_{\textnormal{ns}^+}}(n_{\textnormal{ns}^+}^2)$, where $v_p(m)$ is the exponent of the prime $p$ in the prime factorization of $m$.
\end{proof}

\section{Computing the number of points}\label{sec:basic}

%\subsection{Counting points of modular curves using Hecke eigenvalues} \label{sec:Hecke} 

Let $p$ and $\ell$ be distinct primes. Let $C$ be a smooth, projective, algebraic curve over $\QQ$  with good reduction over $\FF_p$. Let $\mathrm{Ta}_\ell(\mathrm{Jac}(C))$ be the Tate module of $C$, i.e., the inverse limit of the $\ell^m$-torsion group $\mathrm{Jac}(C)[\ell^m]$, and let $V_\ell=\mathrm{Ta}_\ell(\mathrm{Jac}(C))\otimes_{\ZZ_\ell} \QQ_\ell$ be the $\QQ_\ell$-vector space associated to $\mathrm{Ta}_\ell(\mathrm{Jac}(C))$. The trace of the Frobenius automorphism $\mathrm{Frob}_{p^k}$ acting on $V_\ell$ satisfies:
\begin{equation}\label{eq:numpts}
\mathrm{tr}(\mathrm{Frob}_{p^k}|V_\ell)= p^k+1-\#C(\mathbb{F}_{p^k}),
\end{equation}
for every positive integer $k$.
(see \cite[Theorem 11.1]{Mil86jac}). When $C$ is a modular curve associated to a group of matrices containing the scalar matrices, the Eichler-Shimura relation (see \cite[Theorem 8.7.2]{DS05}) implies that the trace of Frobenius can be obtained from a linear combination of traces of Hecke operators, as we explain in \Cref{sec:asympt}.

We are interested in applying \Cref{eq:numpts} when $C=X(n_0$, $n_{0^+}, n_{\textnormal{s}},n_{\textnormal{s}^+}, n_{\textnormal{ns}}, n_{\textnormal{ns}^+})$ as in \Cref{eq:modular-curve-symbol} and $p$ does not divide $n=n_0$ $n_{0^+} n_{\textnormal{s}}n_{\textnormal{s}^+} n_{\textnormal{ns}} n_{\textnormal{ns}^+}$. Using \Cref{th:mixedjac}, the number of points over $\FF_p$ can be computed in terms of the eigenvalues of the Hecke operator $T_p\in \End(J_0^{\textnormal{new}}(d))$, whose level $d$ divides $n$. By \cite[Theorem 6.6.6]{DS05}, we have that $J_0^{\textnormal{new}}(d)$ is isogenous to a direct sum $\bigoplus_f A_f^{\sigma_0(d/d')}$ of Abelian varieties $A_f$ associated to the Galois orbits of the normalized eigenforms $f\in \mathcal S_2(\Gamma_0(d'))$, where $d'\mid d$. In the following part of this section
we explain how to compute $\mathrm{tr}(\mathrm{Frob}_p|V_\ell)$ explicitly when $C=X_H$. This, by \Cref{eq:numpts} above, allows to compute the number of points of $X_H$ on $\mathbb{F}_p$. Finally, we explain how to compute the number of points of $X_H$ on every finite field (not only prime fields). The following lemma is well known, but we report the proof for convenience of the reader.
\begin{lem}\label{snowden}
Let $n$ be a positive integer, let $p$, $\ell$ be distinct primes not dividing $n$ and let
$V_\ell=\mathrm{Ta}_\ell(J_0(n)_{\FF_p})\otimes_{\ZZ_\ell} \QQ_\ell$ be the $\QQ_\ell$-vector space associated to the Tate module of $X_0(n)_{\FF_p}$. Then
\begin{equation*}
\mathrm{tr}_{\mathbb T}(\mathrm{Frob}_p|V_\ell)= T_p,
\end{equation*}
where $T_p$ is the Hecke operator associated to the prime $p$ and the trace is taken on $V_\ell$ as $\mathbb{T}_{\QQ_\ell}$-module with $\mathbb{T}_{\QQ_\ell}:=\mathbb{T}_\ZZ\otimes_\ZZ \QQ_\ell$ and $\mathbb{T}_\ZZ$ denotes the Hecke algebra $\ZZ[T_n: n>0]$.
\end{lem}
\begin{proof}
We have that $V_\ell$ is a free rank 2 module over $\mathbb{T}_{\QQ_\ell}$ (see \cite[Lemma 9.5.3]{DS05}). Consider the Weil pairing $\langle\cdot,\cdot\rangle\colon V_\ell\times V_\ell \rightarrow \QQ_\ell(1)$, see \cite[before Lemma 16.2, p. 132]{Mil86abel} where we identify $J_0(n)$ with its dual variety via the principal polarization. Let $V_\ell^*$ be the dual module of $V_\ell$  where for every linear operator $T\colon V_\ell\to V_\ell$ we have the dual action on $V_\ell^*$, i.e., $T\varphi:=\varphi\circ T$, for every $\varphi\in V_\ell^*$. Let $\Phi\colon V_\ell\rightarrow V_\ell^*$ be the isomorphism sending $x$ in the map $\Phi_x\colon V_\ell\rightarrow \QQ_\ell$ defined by $\Phi_x(y):=\langle y,x\rangle$. Each Hecke operator $T_n$ is self-adjoint with respect to the Weil pairing, see \Cref{lem:T_adj} below. This implies that $\Phi$ commutes with the natural action of $\mathbb{T}_{\QQ_\ell}$ on each side, i.e., $\Phi$ is a isomorphism of $\mathbb{T}_{\QQ_\ell}$-modules.

The Verschiebung $\mathrm{Ver}_p=[p]\mathrm{Frob}_p^{-1}$ is the dual isogeny of $\mathrm{Frob}_p$ on Jacobians and it is the adjoint of $\mathrm{Frob}_p$ with respect to Weil pairing, in fact: $\langle x,\mathrm{Ver}_p(y)\rangle=\langle x,[p]\mathrm{Frob}_p^{-1}(y)\rangle=\langle x,\mathrm{Frob}_p^{-1}(y)\rangle^p=\mathrm{Frob}_p(\langle x,\mathrm{Frob}_p^{-1}(y)\rangle)=\langle \mathrm{Frob}_p(x),\mathrm{Frob}_p(\mathrm{Frob}_p^{-1}(y))\rangle=\langle\mathrm{Frob}_p(x),y\rangle$. Hence $\Phi\circ \mathrm{Frob}_p=\mathrm{Ver}_p\circ\Phi$. This implies that $\mathrm{tr}_{\mathbb T}(\mathrm{Frob}_p|V_\ell) = \mathrm{tr}_{\mathbb T}(\mathrm{Ver}_p|V_\ell^*)$. But, with respect to a fixed basis on $V_\ell$ and its dual basis on $V_\ell^*$, the matrix of $\mathrm{Ver}_p$ on $V_\ell^*$ is the transpose of the matrix of $\mathrm{Ver}_p$ on $V_\ell$, so $\mathrm{tr}_{\mathbb T}(\mathrm{Frob}_p|V_\ell)=\mathrm{tr}_{\mathbb T}(\mathrm{Ver}_p|V_\ell)$. Applying the trace to the Eichler-Shimura relation $\mathrm{Frob}_p+\mathrm{Ver}_p=T_p$ (see \cite[Theorem 8.7.2]{DS05}), we get $2\mathrm{tr}_{\mathbb T}(\mathrm{Frob}_p|V_\ell)=\mathrm{tr}_{\mathbb T}(T_p|V_\ell)=2T_p$ and, hence $\mathrm{tr}_{\mathbb T}(\mathrm{Frob}_p|V_\ell)=T_p$.
\end{proof}

\begin{lem}\label{lem:T_adj}
For all pairwise coprime positive integers $N,n,\ell$, with $\ell$ prime, the Hecke operator $T_n$ is self-adjoint with respect to the Weil pairing on $J_0(N)[\ell^m]$.
\end{lem}
\begin{proof}
Since $T_n$ can be written as a polynomial in the operators $T_p$ for primes $p$ dividing $n$ and all these $T_p$'s commute with each other, then it is enough to treat the case where $n=p$ is prime. 

We focus on the action of $T_p$ on divisors supported on non-cuspidal points. Indeed, since the number of cusps is finite, each element of $J_0(N)$ can be represented as a divisor supported in the non-cuspidal locus. We recall that the non-cuspidal points of $X_0(N)$ parametrize data $(E,G)$, where $E$ is an elliptic curve and $G$ is a cyclic subgroup of order $N$ of $E[N]$ and the non-cuspidal points of $X_0(Np)$ parametrize data $(E,\varphi\colon E\rightarrow E', G)$, for $E,G$ as before and $\varphi$ an isogeny of degree $p$. With this notation, let $p_1,p_2\colon X_0(Np)\rightarrow X_0(N)$ be the maps defined by 
$$
p_1(E,\varphi\colon E\rightarrow E',G)=(E,G),\quad  p_2(E,\varphi\colon E\rightarrow E',G)=(E',\varphi(G)),
$$ 
and extended by continuity on the whole curve. Then $T_p = (p_2)_*p_1^*$ as an endomorphism of the Jacobian. 

We recall that for each non-constant map of curves $f\colon C \to C'$, the maps $f^*$ and $f_*$ between the two Jacobians are adjoint with respect to the Weil pairings. Indeed, given divisors $D,D'$ respectively on $C,C'$ that define $N$-torsion points in the respective Jacobians and such that $f_*(D)$ is disjoint from $D'$, take rational functions $g,g'$ such with $\mathrm{div}(g) = ND$ and $\mathrm{div}(g') = N D'$, so that 
$$ \langle D, f^*D'\rangle = \frac{g(f^*D')}{(g'\circ f)(D)} = \frac{ \mathrm{(Norm}_f g)(D')}{(g'( f_*(D))}  =   \langle f_* D, D'\rangle, $$
where $\langle, \rangle$ denotes the Weil pairing on the $N$-torsion subgroup and $\Norm_f$ is the norm of the finite extension of function fields $f^\#\colon  K(C') \hookrightarrow K(C)$, namely for each rational function $h$ on $C'$ we denote $\Norm_f(h) = \det_{K(C)}(\cdot h\colon K(C') \to K(C'))\in K(C) $.

As a consequence of the above general fact, the adjoint of $T_p$ is $T_p^*=(p_1)_*p_2^*$, i.e., it is induced by the transposed correspondence of $T_p$.
Hence, to prove our statement, it is enough proving that the image of the map $(p_1,p_2)\colon  X_0(Np)\rightarrow X_0(N)\times X_0(N)$ is symmetric. Indeed, given $x=(E,G)$, $y=(E',G')$ in $X_0(N)\times X_0(N)$, the point $(x,y)$ lies in the image of $(p_1,p_2)$ if and only if there exists an isogeny $\varphi\colon E\to E'$ of degree $p$ such that $\varphi(G)=G'$. If this happens, then the dual isogeny
$\hat{\varphi}\colon E'\to E$  sends $G'$ into $G$, since $G=[p]G=\hat{\varphi}(\varphi(G))=\hat{\varphi}(G')$. In other words whenever $(x,y)$ lies in the image of $(p_1,p_2)$, then $(y,x)$ lies in the same image, i.e., the image of $(p_1,p_2)$ is symmetric.
\end{proof}

The following proposition explain how to compute $\mathrm{tr}(\mathrm{Frob}_p|V_\ell)$ explicitly when the abelian variety considered is not necessarily a Jacobian but it is a product of $A_f$ for $f\in\mathcal S_2^{\textnormal{new}}(\Gamma_0(n_f))$, where $n_f\in\ZZ_{>0}$ is the level of $f$. It implies the case we are interested in: $C=X_H$ for $H$ as described in the statement of \Cref{th:mixedjac}.

\begin{prop}\label{prop:eigenvalues_frob_Tp}
Let $\mathcal F$ be a finite subset of the set of the Galois orbits of normalized eigenforms of $\bigcup_{n>0}\mathcal S_2^{\textnormal{new}}(\Gamma_0(n))$, and let $J:=\prod_{[f]\in\mathcal F}A_f^{m_f}$, where $A_f$ is the  abelian variety associated to Galois orbit of $f$ (see \cite[Definition 6.6.3]{DS05}) and $m_f\in\ZZ_{>0}$. Moreover, let $V_\ell=\mathrm{Ta}_\ell(J)\otimes_{\ZZ_\ell} \QQ_\ell$ be the $\QQ_\ell$-vector space associated to the Tate module of $J$. 
Then the characteristic polynomial of $\mathrm{Frob}_p$ acting on $V_\ell$ is 
$$
\prod_{[f]\in\mathcal F}\prod_{h\in [f]} (x^2-a_p(h)x+p)^{m_f}\,,
$$
where $a_p(h)$ is the $p$-th Fourier coefficient of $h$.
\end{prop}
\begin{proof}
It is enough to treat the case  $\mathcal F = \{[f]\}$ for $f$ an eigenform of level $n$. Let
$$
S := \bigoplus_{h\in[f]} \CC  h \quad  \subset \mathcal S_2^\new(\Gamma_0(n)),
$$
and let $\mathbb T$ be the Hecke algebra of $A_f$ defined as the maximal quotient of the Hecke algebra $\mathbb T'$ of $J_0(n)$ such that the action of $\mathbb T'$ on $A_f$ acts through $\mathbb T$. Also, for each ring $R$, let $\mathbb T_R:= \mathbb T\otimes R$. 
%follows: let $\mathbb T'$ be the Hecke algebra of $J_0(n)$ and let $I_f$ be the ideal of operators $T\in\mathbb T'$ such that $T(f) = 0$, hence $\mathbb T := \mathbb T'/I_f$, so that $\mathbb T$ be the smallest quotient of $\mathbb T'$ such that the action of $\mathbb T'$ on $A_f$ factors through $\mathbb T$. 

By \cite[Lemma 9.5.3]{DS05}, $V_\ell$ is a free module of rank $2$ over $\mathbb T_{\QQ_\ell}$, and the characteristic polynomial of $\Frob_p$ as a  
$\mathbb T_{\QQ_\ell}$-linear endomorphism is $x^2-T_px+p$. Indeed, by \Cref{snowden} 
this characteristic polynomial is of the form $x^2-T_px+a$ for a certain 
$a\in \mathbb T_{\QQ_\ell}$. Since $\mathbb T_{\QQ_\ell}$ acts faithfully on
$V_\ell$, the element $a$ is the unique element of $\mathbb T_{\QQ_\ell}$ such that 
$\Frob_p^2-T_p\Frob_p+a  = 0$, and $a=p$ satisfies this property because, by 
Eichler-Shimura Relation and denoting by $\mathrm{Ver}_p$ the Verschiebung operator, we have
$$
\begin{aligned}
\Frob_p^2 - T_p \Frob_p + p &= \Frob_p^2 - (\Frob_p+\mathrm{Ver}_p)\Frob_p + \Frob_p\mathrm{Ver}_p = \\
&=(\Frob_p - \mathrm{Ver}_p)(\Frob_p-\Frob_p) = 0.
\end{aligned}
$$
In particular, the characteristic polynomial of $\Frob_p$ as a $\QQ_\ell$-linear endomorphism of $V_\ell$ is
\[
\begin{aligned}
& \det {}_{\QQ_\ell[x]}(\Frob_p-x| V_\ell \otimes_{\QQ_\ell}\QQ_\ell[x]) = \Norm_{\mathbb{T}_{\QQ_\ell}[x]/\QQ_\ell[x]} \Big( \det{}_{\mathbb{T}_{\QQ_\ell}[x]}(\Frob_p-x| V_\ell \otimes_{\QQ_\ell}\QQ_\ell[x])\Big) =
\\ & \quad =  
\Norm_{\mathbb{T}_{\QQ_\ell}[x]/\QQ_\ell[x]} (x^2 - T_px +p) =  
\Norm_{\mathbb{T}_{\QQ}[x]/\QQ[x]} (x^2 - T_px +p) = \\
& \quad = \Norm_{\mathbb{T}_{\CC}[x]/\CC[x]} (x^2 - T_px +p),
\end{aligned}
\]
where, given a ring extension $A\subset B$ such that $B$ is a finite, free $A$-module, we denote $\Norm_{B/A}(t) := \det_A(\cdot t\colon B \to B)$, i.e., the determinant of the multiplication by $t$ inside $B$ as a $A$-module endomorphism. Since $S$ is a free $\mathbb T_{\CC}$-module of rank $1$, then $S\cong \mathbb T_\CC$ as $\mathbb T_\CC$-modules, hence,  
$$
\Norm_{\mathbb{T}_{\CC}[x]/\CC[x]} (x^2 - T_px +p) = \det{}_{\CC[x]}(x^2 - T_px +p | S\otimes_\CC \CC[x]) = \prod_{h\in [f]} (x^2-a_p(h)x+p),
$$
where the last equality is true because the elements $h\in [f]$ are a basis of eigenvectors for $T_p$ in $S$, with eigenvalue $a_p(h)$.
\end{proof}

%\subsection{An algorithm for $|X(n_0,n_{0^+}, n_{\textnormal{ns}}, n_{\textnormal{ns}^+})(\mathbb{F}_q)|$} \label{sec:alg} 

Using Theorem \ref{th:mixedjac} and Formula \ref{snowden} above we can compute $\#X(n_0,n_{0^+}, n_{\textnormal{ns}}, n_{\textnormal{ns}^+})(\mathbb{F}_q)$ for all prime power $q=p^k$ not dividing $n=n_0 n_{0^+} n_{\textnormal{ns}} n_{\textnormal{ns}^+}$. We concentrate on these cases because of the isomorphism of \Cref{rem:s_non_serve}. The algorithm runs as follow: 
\begin{Algorithm} \emph{Computing the number of points of $X(n_0,n_{0^+},n_\ns,n_{\ns^+})$ over $\mathbb F_q$.}\label{alg:alg}
\begin{enumerate}
    \item For every $d|n$, compute $d_0:=\gcd(n_0,d)$ and similarly $d_0^+$, $d_{\textnormal{ns}}$ and $d_{\textnormal{ns}^+}$ and let $D:=d_0d_0^+d_{\textnormal{ns}}^2d_{\textnormal{ns}^+}^2$.
    \item For every $d|n$, choose a basis $\mathcal B(d)$ of new eigenforms of level $D$ and let $\mathcal B:=\bigcup_{d\mid n} \mathcal B(d)$ be the union of these bases.
    \item For every $d\mid n$ and for each prime $\ell\mid d_{\textnormal{ns}}^+$, compute the eigenvalue of $f\in\mathcal B(d)$ %of level the maximal power $\ell^r$ dividing $d_{\textnormal{ns}}^+$,
    with respect to the Atkin-Lehner operator $w_\ell$ and discard from $\mathcal B$ the $f$'s with eigenvalue $-1$ (these data are available on the database \cite{lmfdb} for $D\leq 10000$).
    \item Compute the Hecke eigenvalue $a_p(f)$ for every $f\in\mathcal B$ (these data are available on the database \cite{lmfdb} for $D\leq 10000$).
    \item For each $f\in\mathcal B$, compute $m_f:=\sigma_0(n_0/d_0)\tilde\sigma_0(n_{0^+}/d_{0^+},f)$, where $\sigma_0(m)$ is the number of positive divisors of $m$ and $m\mapsto \tilde\sigma_0(m,f)$ is the multiplicative function defined by
    \[
\tilde\sigma_0(p^r,f):=\begin{cases}
\frac{1}{2}(r+1), & \text{if }r\text{ is odd}, \\
\frac{r}{2}+1, & \text{if }r\text{ is even and } w_pf=f, \\
\frac{r}{2}, & \text{if }r\text{ is even and } w_pf=-f.
\end{cases}
\]
    \item If $k=1$, then compute Frobenius traces acting on $A_f$: $\mathrm{tr}(\mathrm{Frob}_p|A_f)=\sum_{h\in [f]} a_p(h)$, where $[f]$ is the Galois orbit of $f$.
    \item If $k>1$, then compute the (complex) roots $\alpha(h)$ and $\beta(h)$ of the polynomial $x^2-a_p(h)x+p=0$, for each $h\in [f]$. Then compute $\mathrm{tr}(\mathrm{Frob}_q|A_f)=\sum_{h\in [f]}(\alpha(h)^k+\beta(h)^k)$, using \Cref{eq:numpts}
    \item Finally, $\#X(n_0,n_{0^+}, n_{\textnormal{ns}}, n_{\textnormal{ns}^+})(\FF_q)=q+1-\sum_{[f]} m_f \mathrm{tr}(\mathrm{Frob}_q|A_f)$, where the sum is taken over the distinct orbits for $f\in \mathcal{B}$.
\end{enumerate}
\end{Algorithm}

\begin{comment}
\begin{ese}\label{livello23}
As an example we apply our algorithm for $q=3^2$ and $n_0=23$. 
The genus of $X_0(23)=X(23,1,1,1)$ is $2$, in fact this is the dimension of $\mathcal{S}_2(\Gamma_0(23))$. 
A basis of this space is given by $f=q + a\cdot q^2 + (-2a - 1)q^3 + (-a - 1)q^4 + 2a\cdot q^5 + (a - 2)q^6 + (2a +2)q^7 +\ldots$ and its Galois conjugated $\bar{f}$, where $a$ is a root of $x^2+x-1$. 
The characteristic polynomial of the Hecke operator $T_3$ acting on the abelian variety $A_f$ associated to $f$ is 
%$P(x)=x^4 - 10x^2+ 25$. $P(x)=(x^2-5)^2$. 
$P(x)=x^2-5$. The roots are $t_1=\sqrt{5}$ and $t_2=-\sqrt{5}$. %have multiplicity 2 so we take them once so the genus of $X_0(23)$ is $2$. 
Consider the polynomial $(x^2-t_1x+p)(x^2-t_2x+p)$. Its roots are
$\alpha_{1,2}=\frac{\sqrt{5}}{2} \pm \frac{\sqrt{7}}{2}i$ and $\alpha_{3,4}=-\frac{\sqrt{5}}{2} \pm \frac{\sqrt{7}}{2}i$. It follows that $|X_0(23)(\FF_{9})|=9+1-(\alpha_1^2+\alpha_2^2+\alpha_3^2+\alpha_4^2)=12$.
\end{ese}
\end{comment}
\begin{ese}\label{livello163+}
We compute the number of points of $X(1,163,1,1,1,1)=X_0^+(163)$ over $\FF_{2^2}$ without using the explicit equation in \cite{Mer18}. Consider the space of cusp forms $\mathcal{S}_2^\new(\Gamma_0(163))$. This space has dimension 13 and let $f_1,f_2,f_3$ be representatives under the Galois action. The class of $f_1$ has just one element, the class of $f_2$ has 5 elements and the class of $f_3$ has 7 elements. Since $f_3$ has negative eigenvalue, its class is discarded.
We compute $m_{f_1}=m_{f_2}=1$, hence the genus of $X_0^+(163)$ is $6$. 
%The characteristic polynomials of the Hecke operator $T_2$ acting on the abelian variety $A_f$ associated to $f$ for $f=f_1$ and $f=f_2$ are, respectively, $P_1(x)=x$ and %$P_1(x)=x^2$ and $P_2(x)=x^{10} + 10x^9 + 31x^8 - 173x^6 - 244x^5 + 159x^4 + 498x^3 + 166x^2 - 96x+ 9$. $P_2(x)=(x^5 + 5x^4 + 3x^3 - 15x^2 - 16x + 3)^2$. 
%$P_2(x)=x^5 + 5x^4 + 3x^3 - 15x^2 - 16x + 3$.
%The only root of $P_1$ is 
The Hecke eigenvalues are $a_2(f_1)=0$
and $a_2(h)\approx -2.711, -2.484,-1.634,0.163,1.666$, for $h\in [f_2]$. %all with multiplicity 2, so we take them once. 
%As in Example \ref{livello23} 
We consider the polynomial $(x^2-a_2(f_1)x+2)\prod_{h\in[f_2]}(x^2-a_2(h)x+2)$. %for $j=1,\ldots,6$. 
We denote its roots by $\alpha_j$ and $\beta_{j}$, where
$\alpha_{1}\approx 1.414i$,
$\alpha_{2}\approx -1.355+ 0.402i$,
$\alpha_{3}\approx -1.242 + 0.675i$,
$\alpha_{4}\approx -0.817+ 1.154i$,
$\alpha_{5}\approx 0.081 + 1.411i$, 
$\alpha_{6}\approx 0.833+ 1.1425i$ and $\beta_{j}=\bar \alpha_{j}$ is the complex conjugated of $\alpha_j$. 
Hence $\#X_0^+(163)(\FF_{2^2})=4+1-\sum_{j=1}^6 (\alpha^2_j+\beta^2_{j})=4+1+5=10$.
\end{ese}

\begin{ese}\label{livello_17,3}
We compute the number of points of $X(1,17,1,1,1,3)$ over $\FF_2$. 
We have to consider the space of cusp forms $\mathcal{S}_2^\new(\Gamma_0(D))$ at levels $D=1,9,17,153$. 
The first two spaces are trivial. The third one has dimension 1, but the corresponding representative under the Galois action must be discarded, since the eigenvalue of Atkin-Lehner operator $w_{17}$ is $-1$.
Finally, the last space has 5 representatives under the Galois action and we denote them by $f_1,f_2,f_3,f_4,f_5$. Four of these classes should be discarded since they have at least a negative eigenvalue.
We denote by $f$ the unique representative whose both eigenvalues of $w_3$ and $w_{17}$ are positive. 
Then $A_f$ has dimension 1 and $m_f=\sigma_0(1)\sigma^*_0(1)=1$, so it follows that the genus of $X(1,17,1,1,1,3)$ must be 1. 
%The characteristic polynomial of the Hecke operator $T_2$ acting on the abelian variety $A_f$ associated to $f$ is $P(x)=x+2$. 
The eigenvalue of $T_2$ acting on $f$ is $a_2(f)=-2$.
%Denote by $t$ the unique root of $P(x)$, so we consider the polynomial $x^2-tx+2=x^2+2x+2$.
%We denote the roots of  by $\alpha_1$ and $\alpha_2$, then $\alpha_1+\alpha_2=-2$ and 
Then $\#X(1,17,1,1,1,3)(\FF_{2})=2+1+2=5$. 
This is the maximum possible value for an elliptic curve over $\FF_2$ according to \cite{manypoints}.
\end{ese}

\begin{ese}\label{livello_5,3,2}
We compute the number of points of $X(5,3,1,1,2,1)$ over $\FF_{7^2}$. We have to consider the space of cusp forms $\mathcal{S}_2^\new(\Gamma_0(D))$ at levels $D=1,3,4,5,12,15,20,60$. Only the spaces corresponding to $D=15$ and $D=20$ are non-trivial, more precisely, $\mathcal{S}_2^{\new}(\Gamma_0(D))$ has dimension 1 in both these cases. We denote by $f_1$ and $f_2$, respectively, the corresponding representative under the Galois action. In both cases, the action under the Atkin-Lehner involution $w_3$ gives positive eigenvalues so both  $f_1$ and $f_2$ must be considered.
Since $m_{f_1}=\sigma_0(1)\sigma^*_0(1)=1$ and $m_{f_2}=\sigma_0(1)\sigma^*_0(3)=1$, it follows that the genus of $X(5,3,1,1,2,1)$ must be 2. 
%The characteristic polynomials of the Hecke operator $T_7$ acting on the abelian varieties $A_{f_1}$ and $A_{f_2}$ are $P_1(x)=x$ and $P_2(x)=x-2$, respectively. 
The eigenvalue of $T_7$ acting on $f_1$ and $f_2$ are $0$ and $2$, respectively.
We denote by $\alpha_1,\beta_1,\alpha_2,\beta_2$ the roots of the polynomial $(x^2+7)(x^2-2x+7)$. Then 
$\alpha_{1}\approx 2.646i$,
$\beta_{1}\approx -2.646i$,
$\alpha_{2}\approx 1 + 2.499i$,
$\beta_{2}\approx 1 - 2.499i$
and $\#X(5,3,1,1,2,1)(\FF_{7^2})=7^2+1- \sum_{i=1}^2(\alpha_i^2+\beta_i^2)= 49+1+24=74$. 
\end{ese}

\section{Asymptotic estimates}\label{sec:asympt}

In this section we estimate the number of points on the reductions over prime finite fields of the curves $X = X(n_0,n_{0^+}, n_{\textnormal{s}}, n_{\textnormal{s}^+}, n_{\textnormal{ns}}, n_{\textnormal{ns}^+})$. We start by recalling a classical formula which we then combine with estimates on the trace of Hecke operators. Notice that by \Cref{th:mixedjac}, the space of regular differentials on $X$ can be identified with a direct sum of spaces of newforms of shape $\mathcal S_2^\new(\Gamma_0(d))^W$ for $W$ a certain group of Atkin-Lehner operators and $d$ a suitable positive integer. In particular, this decomposition defines an action of Hecke operators on $\Omega^1_{X/\CC}$.  

\begin{prop}\label{prop:trace_points}
	Given $X = X(n_0,n_{0^+}, n_{\textnormal{s}}, n_{\textnormal{s}^+}, n_{\textnormal{ns}}, n_{\textnormal{ns}^+})$ and defining $T_{p^{-1}}:=0$, 
	we have
	$$
	\# X(\FF_{p^e}) =   
	p^e+1 - \mathrm{tr} (T_{p^e}|\Omega^1_{X/\CC}) + p \, \mathrm{tr}(T_{p^{e-2}}|\Omega^1_{X/\CC}), 
	$$
	for each integer $e\ge 1$ and for every prime $p$ not dividing  $n= n_0 n_{0^+} n_{\textnormal{s}}n_{\textnormal{s}^+} n_{\textnormal{ns}} n_{\textnormal{ns}^+}$.
\end{prop}

\begin{proof}
	Let $\alpha_1, \ldots, \alpha_{g},\beta_1,\ldots,\beta_g$ be the eigenvalues of the Frobenius acting on the Tate module of the Jacobian of $X_{\FF_p}$. By \Cref{prop:eigenvalues_frob_Tp}, up to reordering, we can suppose that $\alpha_i+\beta_{i} = a_p(f_i)$, where $f_1, \ldots, f_g$ are the normalized eigenforms giving a basis of $\Omega^1_{X/\CC}$. Hence  because of \Cref{eq:numpts} and \Cref{prop:eigenvalues_frob_Tp}, it is enough to prove that
	\begin{equation}\label{eq:intermediate_trace}
		\alpha_i^e + \beta_i^e = a_{p^e}(f_i) - pa_{p^{e-2}}(f_i).
	\end{equation}
	The case $e=1$ is true defining $a_{p^{-1}}:=0$. Since $\alpha_i\beta_i = p$ (by \Cref{prop:eigenvalues_frob_Tp}), and $T_{p^2} = T_p^2 - p$ (by \cite[Equation (5.10)]{DS05}), we have
	\[
	a_{p^2}(f_i) -p = a_p(f_i)^2 - 2p = (\alpha_i + \beta_i)^2 - 2\alpha_i\beta_i = \alpha_i^2 + \beta_i^2, 
	\]
	which proves the case $e=2$. For the inductive step we recall that on $\Omega^1_{X/\CC}$ we have
	\[
	T_{p^e} - p T_{p^{e-2}} = %T_p (T_p^{e-1}-p \langle p\rangle T_p^{e-3}) - p (T_{p^{e-2}}-p\langle p\rangle T_{p^{e-4}})  = 
	T_p (T_{p^{e-1}}-pT_{p^{e-3}}) - p (T_{p^{e-2}}-pT_{p^{e-4}}).
	\]
	Hence, supposing that \Cref{eq:intermediate_trace} holds for $e{-}1$ and $e{-}2$, we get 
	\[
	\begin{aligned}
		a_{p^e}(f_i) - p a_{p^{e-2}}(f_i) & = a_p(f_i) (a_{p^{e-1}}(f_i)-pa_{p^{e-3}}(f_i)) - p (a_{p^{e-2}}(f_i)-pa_{p^{e-4}}(f_i)) \\
		& = (\alpha_i + \beta_i)(\alpha_i^{e-1} + \beta_i^{e-1}) - (\alpha_{i}\beta_i)(\alpha_i^{e-2} + \beta_i^{e-2}) = \alpha_i^e + \beta_i^e.
	\end{aligned}
	\]
\end{proof} 
To estimate the traces of Hecke operators on $\mathcal S_2^\new(\Gamma_0(d))^W$ for $W$ a group of Atkin-Lehner operators, %and $d$ a suitable positive integer, 
it is enough to estimate $\mathrm{tr}(T_p w_m | \mathcal S_2^\new(\Gamma_0(d)))$ for each $w_m\in W$. A small remark on notation: given the prime factorization $d= p_1^{e_2}\cdots p_r^{e_r}$, here and in  \Cref{sec:jac},
we sometimes write $w_{p_i}$ meaning $w_{p_i^{e_i}}$, while other times we use the usual notation $w_m$ for $m\mid d$ such that $\gcd(m,d/m)=1$.
%, namely a product of $p_i^{e_i}$'s. 
To estimate the trace of $T_k w_m$, for such a divisor $m$ of $d$, we look at \cite[Proposition~2.8]{Bru}, which states:
\begin{equation}\label{eq:Bru_tr}
	|\mathrm{tr}(T_k w_m| \mathcal S_2^\new(\Gamma_0(d))) -  \delta(\sqrt k\in \ZZ) F(d,m)| < c_0 (k + \sqrt{kd}) \sigma_0^3(d) \sigma_0(k) \log^2(4kd),     
\end{equation}
where $c_0$ is an absolute constant, $\sigma_0$ is the function counting the number of divisors of a natural number, $\delta( \text{condition})$ is $1$ if the condition is true and 0 otherwise, and
$$
F(d,m) := \mu(\sqrt m) \frac{d}{m} \prod_{\ell \in \mathcal{P}( \frac dm)}\Big(1-\frac 1\ell - \frac{\delta(\ell^2|\tfrac{d}{m})}{\ell^2} + \frac{\delta(\ell^3|\tfrac{d}{m})}{\ell^3} \Big),
$$
with $\mathcal P(N)$ the set of prime divisors of an integer $N$, and $\mu$ the M\"oebius function extended so that $\mu(\sqrt m)=0$ if $\sqrt m$ is not an integer. Combining \Cref{eq:Bru_tr} and \Cref{prop:trace_points} we get the following proposition.
\begin{prop}
	There exists an absolute constant $C$ such that, for each curve $X= X(n_0,n_{0^+}, n_{\textnormal{s}}, n_{\textnormal{s}^+}, n_{\textnormal{ns}}, n_{\textnormal{ns}^+})$ of genus $g$ and for each prime power $q=p^e$ coprime to $n_0n_{0^+} n_{\textnormal{s}} n_{\textnormal{s}^+} n_{\textnormal{ns}} n_{\textnormal{ns}^+}$, we have 
	$$
	\begin{cases} 
		\#X(\FF_{q}) <  C (q+\sqrt{gq})g^{\frac{10}{\log\log g}}p\log^3 q, & \text{ if $q$ is not a square,} \\
		%C\big(qg^{\frac{2}{\log\log g}}+ q^{\frac12}g^{\frac12+\frac{12}{\log\log g}}\log q \big) p\log q & \text{ if $q$ is not a square,}\\
		|\#X(\FF_{q}) - (p-1)g | <  C (q+\sqrt{gq})g^{\frac{10}{\log\log g}}p\log^3 q ,
		% C\big(qg^{\frac{2}{\log\log g}}+ q^{\frac12}g^{\frac12+\frac{12}{\log\log g}}\log q \big) p\log q
		 & \text{ if $q$ is a square.}
	\end{cases}
	$$
\end{prop}
\begin{proof}
	We start by slightly simplifying \Cref{eq:Bru_tr}: %when $d<11$ the space $S_2^\new(\Gamma_0(d))$ is trivial, while, for $d\ge 11$, 
	using \cite[Theorem 1]{Rob83sigma0} to bound $\sigma_0$, we get
	\begin{equation}\label{eq:Bru_g}
		|\mathrm{tr}(w_m| \mathcal S_2^\new(\Gamma_0(d))) -  F(d,m)| < 2 c_0 d^{\frac 12 + \frac{4.4 }{\log\log d}} %< c_\epsilon N^{\frac 12 +\epsilon},    
	\end{equation}
	when specializing to $k=1$ in \Cref{eq:Bru_tr}; while for $k=r$ a prime power we have
	\begin{equation}\label{eq:Bru_r}
		|\mathrm{tr}(T_r w_m| \mathcal S_2^\new(\Gamma_0(d))) -  \delta(\sqrt r\in \ZZ) F(d,m)| < 20 c_0  d^{\frac{4.4}{\log\log d}} (r+ \sqrt{rd})\log^3 r.%< c_\epsilon (k^{1 +\epsilon} + (k N)^{\frac 12 +\epsilon}),
	\end{equation}
	Since we are are only interested in the cases where $\mathcal S_2^\new(\Gamma_0(d))\ne 0$, we restrict to the cases $d\ge 11$, hence $\log\log d$ makes sense and is positive. 
	
	We recall that, for each finite group $G$ acting on a finite vector space $V$, the linear operator $\sum_{g\in G} g\colon V \to V$ has trace equal to $\# G \dim (V^G)$. Hence, applying \Cref{eq:Bru_g}, for each group of the form $W = \langle w_{\ell}: \ell \in \mathcal Q\rangle $ with $\mathcal Q$ a subset of the prime divisors of $d$, we can write
	\begin{equation} \label{eq:J0new_approx}
		\dim(J_0(d)^W) = \frac{1}{\# W}\sum_{w_m\in W} \mathrm{tr}(w_m| \mathcal S_2^\new(\Gamma_0(d))) =  F(d,W) + \epsilon(d,W),
	\end{equation}
	where $|\epsilon(d,W)|< 2 c_0 d^{\frac 12 + \frac{4.4 }{\log\log d}}$
	%$$ %|\epsilon_{d,W}| < 2 c_0 d^{\frac 12 + \frac{4.4 }{\log\log d}} %< c_\epsilon N^{\frac 12 +\epsilon}, %$$
	and $F(d,W)$ is defined and can be estimated as follows
	$$
	\begin{aligned}
		& F(d,W)  := \frac{1}{\# W} \sum_{w_m\in W} F(d,m)=\frac{1}{\# W} \sum_{\{\ell_1,\ldots,\ell_r\}\subseteq \mathcal Q} F(d,\ell_1^{e_{\ell_1}}\cdots\ell_r^{e_{\ell_r}})= \\
		& = \frac{d}{\# W}\sum_{\{\ell_1,\ldots,\ell_r\}\subseteq \mathcal Q}  \prod_{\ell\in \{\ell_1,\ldots,\ell_r\}} \frac{\mu(\ell^{e_\ell/2})}{\ell^{e_\ell}} \prod_{\substack{ \ell|d \\ \ell \notin \{\ell_1,\ldots,\ell_r\}}} \Big(1-\frac 1\ell - \frac{\delta(\ell^2|\tfrac{d}{m})}{\ell^2} + \frac{\delta(\ell^3|\tfrac{d}{m})}{\ell^3}  \Big)= \\
		& = \frac{d}{\# W}\prod_{\substack{ \ell|d \\ \ell \notin \mathcal Q}}\Big(1-\frac 1\ell - \frac{\delta(\ell^2|\tfrac{d}{m})}{\ell^2} + \frac{\delta(\ell^3|\tfrac{d}{m})}{\ell^3}  \Big) \cdot \\
        & \cdot \sum_{\{\ell_1,\ldots,\ell_r\}\subseteq \mathcal Q}   \prod_{\ell\in \{\ell_1,\ldots,\ell_r\}} \frac{\mu(\ell^{e_\ell/2})}{\ell^{e_\ell}} \prod_{\substack{ \ell\in \mathcal Q \\ \ell \notin \{\ell_1,\ldots,\ell_r\}}} \Big(1-\frac 1\ell - \frac{\delta(\ell^2|\tfrac{d}{m})}{\ell^2} + \frac{\delta(\ell^3|\tfrac{d}{m})}{\ell^3}  \Big)= \\
		& = \frac{d}{\# W} \prod_{\substack{ \ell|d \\ \ell \notin \mathcal Q}} \Big(1-\frac 1\ell - \frac{\delta(\ell^2|\tfrac{d}{m})}{\ell^2} + \frac{\delta(\ell^3|\tfrac{d}{m})}{\ell^3}  \Big) \prod_{\ell\in \mathcal Q} \Big(1-\frac 1\ell - \frac{\delta(\ell^2|\tfrac{d}{m})}{\ell^2} + \frac{\delta(\ell^3|\tfrac{d}{m})}{\ell^3} + \frac{\mu(\ell^{e_\ell/2})}{\ell^{e_\ell}}\Big)\ge \\
		& \ge \frac{d}{\sigma_0(d)} \prod_{\ell|d } \big(1 - \frac1\ell \big)(1-\frac{1}{\ell^2})^2 \ge \frac{\varphi(d)}{\sigma_0(d)} \zeta(2)^{-2} \ge c_1 d^{1-\frac{1.2}{\log\log d}}, % \ge c_\epsilon d^{1-\epsilon}, 
	\end{aligned}
	$$
	where $\varphi$ is the Euler's totient function, estimated as in \cite[Theorem 327]{HarWri}, and $\zeta$ is the Riemann zeta function. 
	Analogously, \Cref{eq:Bru_r} implies that
	\begin{align} \label{eq:trTq_new_W_apptrox}
		\mathrm{tr}(T_r|\mathcal S_2^\new(\Gamma_0(d))^W) &= \frac{1}{\# W}\sum_{w_m\in W} \mathrm{tr}(T_rw_m|\mathcal S_2^\new(\Gamma_0(d))) = \\ \nonumber
		&=\delta(\sqrt r\in \ZZ) F(d,W) + \epsilon(d,W,r),
	\end{align}
	with 
	$	|\epsilon(d,W,r)| <  20 c_0  d^{\frac{4.4}{\log\log d}} (r+ \sqrt{rd})\log^3 r$. % < c_\epsilon (k + k N)^{\frac 12 +\epsilon}.$
	
	We now look at the curve $X$. Using \Cref{th:mixedjac} to write $\Omega^1_{X/\CC}$ as a sum of spaces of the form $\mathcal S_2^\new(\Gamma_0(d))^W$ and taking linear combinations of  Equations \ref{eq:trTq_new_W_apptrox} and \ref{eq:J0new_approx}, we get
	\begin{equation}\label{eq:g_and_Tr}
		g = F + \epsilon,\qquad \mathrm{tr}(T_r|\Omega^1_{X/\CC}) = \delta(\sqrt r\in \ZZ)F + \epsilon_r,
	\end{equation}
	for
	\[
	\begin{aligned}
		F&= 
		%\sum_{d\mid N} \sigma_0\left(\frac{n_0n_{\textnormal{s}}^2}{d_0d_{\textnormal{s}}}\right)\sigma_0^*\left(\frac{n_{0^+}n_{\textnormal{s}^+}^2}{d_{0^+}d_{\textnormal{s}^+}}\right)
		% F(\langle \prod_{\substack{p|d_{\textnormal{ns}^+}\\p\text{ prime}}}w_p\rangle, d_0d_{0^+}d_{\textnormal{s}} d_{\textnormal{s}^+} d_{\textnormal{ns}}^2 d_{\textnormal{ns}^+}^2)
		\sum_{d} \sigma_0\Big(\frac{n_0n_{\textnormal{s}}^2}{d_0d_{\textnormal{s}}}\Big)
		\sum_{m^2||  \frac{n_{0^+}n_{\textnormal{s}^+}^2}{d_{0^+}d_{\textnormal{s}^+}}}
		\sigma_0^+\Big(\frac{n_{0^+}n_{\textnormal{s}^+}^2}{d_{0^+}d_{\textnormal{s}^+} m^2}\Big)
		F\Big( d_0d_{0^+}d_{\textnormal{s}} d_{\textnormal{s}^+} d_{\textnormal{ns}}^2 d_{\textnormal{ns}^+}^2, \langle w_\ell: \ell|md_{\textnormal{ns}^+} \rangle \Big), \\ 
        \epsilon &= \sum_{d} \sigma_0\Big(\frac{n_0n_{\textnormal{s}}^2}{d_0d_{\textnormal{s}}}\Big)
		\sum_{m^2||\frac{n_{0^+}n_{\textnormal{s}^+}^2}{d_{0^+}d_{\textnormal{s}^+}}}
		\sigma_0^+\Big(\frac{n_{0^+}n_{\textnormal{s}^+}^2}{d_{0^+}d_{\textnormal{s}^+} m^2}\Big)
		\epsilon\Big(d_0d_{0^+}d_{\textnormal{s}} d_{\textnormal{s}^+} d_{\textnormal{ns}}^2 d_{\textnormal{ns}^+}^2, \langle w_\ell : \ell|md_{\textnormal{ns}^+}\rangle \Big), \\ 
        \epsilon_r &= 
		\sum_{d} \sigma_0 \Big(\frac{n_0n_{\textnormal{s}}^2}{d_0d_{\textnormal{s}}}\Big)
		\sum_{m^2||  \frac{n_{0^+}n_{\textnormal{s}^+}^2}{d_{0^+}d_{\textnormal{s}^+}}}
		\sigma_0^+\Big(\frac{n_{0^+}n_{\textnormal{s}^+}^2}{d_{0^+}d_{\textnormal{s}^+} m^2}\Big)
		\epsilon\Big(d_0d_{0^+}d_{\textnormal{s}} d_{\textnormal{s}^+} d_{\textnormal{ns}}^2 d_{\textnormal{ns}^+}^2, \langle w_\ell : \ell|md_{\textnormal{ns}^+} \rangle, r\Big),
	\end{aligned}
	\]
	where the external sums are indexed over $d$ equal to the product of  %= d_0 d_{0^+} d_{\textnormal{s}} d_{\textnormal{s}^+} d_{\textnormal{ns}} d_{\textnormal{ns}^+}$ for certain divisors 
	$d_0, d_{0^+}, d_{\textnormal{s}}, d_{\textnormal{s}^+}, d_{\textnormal{ns}}, d_{\textnormal{ns}^+}$, which vary across the divisors respectively of $n_0, n_{0^+}, n_{\textnormal{s}}^2, n_{\textnormal{s}^+}^2, n_{\textnormal{ns}}, n_{\textnormal{ns}^+}$, while the internal sums of the form $\sum_{m^2||k}$ are indexed over the positive integers $m$ such that $m^2$ divides $k$ and $\gcd(k, \tfrac{k}{m^2})=1$, and finally $\sigma_0^+$ is the multiplicative function such that $\sigma_0^+(\ell^e) = \lceil {\tfrac e2} \rceil$ for each prime power $\ell^e$.
	Under the notation $N = n_0 n_{0^+} n_{\textnormal{s}}^2 n_{\textnormal{s}^+}^2 n_{\textnormal{ns}}^2 n_{\textnormal{ns}^+}^2$, using \cite[Theorem 1]{Rob83sigma0} and the previous estimates on $F(d,W)$, $\epsilon(d,W)$ and $\epsilon(d,W,r)$, we get
	$$
	\begin{aligned}
		F &\ge F(N, \langle w_\ell:\ell| n_{\ns^+} \rangle) \ge c_1 N^{1-\frac{1.2}{\log\log N}} , %\ge c_\epsilon (N)^{1-\epsilon}, 
		\\ 
        |\epsilon | &\le \sigma_0(N)^3 \max\{|\epsilon(d, W)| : d|N, W < \langle w_\ell:\ell|d\rangle \} <  2 c_0 N^{\frac 12 + \frac{7.6}{\log\log N}}, \\ 
        |\epsilon_r | &\le \sigma_0(N)^3 \max\{  \epsilon(d, W, r)  : d|N, W < \langle w_\ell:\ell|d \rangle\} \le 20 c_0 N^{\frac{7.6}{\log\log N}}  (r+ \sqrt{rN})\log^3 r. %  < c_\epsilon (k \sqrt{n''})^{1 +\epsilon}.
	\end{aligned}
	$$
	Plugging the first two bounds above into \Cref{eq:g_and_Tr}, for appropriate constants $c_2$ and $c_3$, we have
	\[
	g> c_2 N^{1-\frac{1.2}{\log\log N}}
	\]
	that implies
	\[
	N< c_3 g^{1+\frac{1.4}{\log\log g}}.
	\]
	Again using \Cref{eq:g_and_Tr} and the bounds on $\epsilon, \epsilon_r$, we get, for a suitable $c_4$, that
	\[
	\begin{aligned}
		|\mathrm{tr}(T_r|\Omega^1_{X/\CC}) - \delta(\sqrt r\in \ZZ)g| & =  | \epsilon_r - \delta(\sqrt r\in \ZZ)\epsilon | \\
		& <  2 c_0 N^{\frac 12 + \frac{7.6}{\log\log N}}  + 20 c_0 N^{\frac{7.6}{\log\log N}}  (r+ \sqrt{rN})\log^3 r \\
		& < c_4 (r+ \sqrt{r}g^{\frac 12+\frac{0.7}{\log\log g}}) g^{\frac{8.4}{\log\log g}}\log^3 r.
	\end{aligned}
	\]
	This estimate, together with \Cref{prop:trace_points}, implies the desired result.
\end{proof}

The above proposition implies that, for a fixed $q$, we expect a large number of points on our curves mostly for $q=p^2$. Anyway, the estimates are only relevant for $N$ large with respect to $q$, and indeed we found records also for fields of the form $p^5$.

%\begin{align}\label{Cartan-Borel-isogeny-fine}
%&\mathrm{Jac}(X_H)/\langle w_{p_{0^+}}, \text{ for $p_{0^+}$ prime dividing }n_{0^+}\rangle \sim \\ \nonumber
%&\sim\prod_{d\mid N}J_0^\new(d_0d_{0^+}d_{\textnormal{s}} d_{\textnormal{s}^+} d_{\textnormal{ns}}^2 d_{\textnormal{ns}^+}^2)^{\sigma_0\left(\frac{n_0n_{\textnormal{s}}^2}{d_0d_{\textnormal{s}}}\right)\sigma_0^*\left(\frac{n_{0^+}n_{\textnormal{s}^+}^2}{d_{0^+}d_{\textnormal{s}^+}}\right)\displaystyle\prod_{\substack{p|d_{\textnormal{ns}^+}\\p\text{ prime}}}w_p},
%\end{align}

\begin{rem}\label{rem:supersingular}
As first shown in \cite{Ogg74} (see also \cite[Section~1]{TVZ82} or \cite[Section~3.3]{BGGP05}),
the case $q=p^2$ is known to provide many rational points in modular curves, thanks to the presence of supersingular points: the $j$-invariant of a supersingular elliptic curve $E$ lies in $\FF_{p^2}$ and, up to twisting, the absolute Galois of $\FF_{p^2}$ acts diagonally on the torsion of $E$.
For example, the bounds in \cite[Lemma 3.20]{BGGP05} imply that, for each congruence subgroup $\Gamma\subset\SLZ$, and for each prime $p$ coprime to the level, we have 
\begin{equation}\label{eq:supersing_gamma}
X_\Gamma(\FF_{p^2})>(p-1)(g_\Gamma-1),
\end{equation} where $g_\Gamma$ is the genus of $X_\Gamma$. This inequality can  easily be extended to our curves. Indeed, let $X = X(n_0,n_{0^+}, n_{\textnormal{s}}, n_{\textnormal{s}^+}, n_{\textnormal{ns}}, n_{\textnormal{ns}^+})$, with genus $g$ and let $p$ be a prime not dividing $n_0n_{0^+} n_{\textnormal{s}} n_{\textnormal{s}^+} n_{\textnormal{ns}} n_{\textnormal{ns}^+}$. Then we can write $X = X_\Gamma/G$ for a suitable choice of $\Gamma$ and of a group of automorphisms $G$. The Riemann-Hurwitz's formula implies that $g$ satisfies $g-1 \le (g_\Gamma-1)/(\# G)$. The inclusion $X(\FF_{p^2}) \supset  X_\Gamma(\FF_{p^2})/G$, together with orbit counting implies that  $X(\FF_{p^2})$ has cardinality at least $\# X_\Gamma(\FF_{p^2})/\#G$. Using also \Cref{eq:supersing_gamma} we get 
$
X(\FF_{p^2}) > (p-1)(g-1).
$
\end{rem}

\section{Greatest Hits}\label{sec:app}

We implemented \Cref{alg:alg} and applied it to all the curves for which we could compute the number of points using the numerical data available at \cite{lmfdb}.
In \Cref{tab:record1} and \Cref{tab:record2} we list all the improved bounds, which we found, for the number of points of $X(n_0,n_{0^+}, n_{\textnormal{ns}}, n_{\textnormal{ns}^+})$ over $\mathbb{F}_{q}$, where $n_0n_{0^+} n_{\textnormal{ns}}^2 n_{\textnormal{ns}^+}^2\le 10000$ and $q=p^k$ is a prime power with $p<20$ and $k\le 5$. 
As explained in \Cref{rem:s_non_serve}, the search along this set of curves gives the same results over the bigger set of 
$X(n_0,n_{0^+}, n_{\textnormal{s}}, n_{\textnormal{s}^+}, n_{\textnormal{ns}}, n_{\textnormal{ns}^+})$ curves since
\[
\#X(n_0,n_{0^+}, n_{\textnormal{s}}, n_{\textnormal{s}^+}, n_{\textnormal{ns}}, n_{\textnormal{ns}^+})(\FF_q)=\#X(n_0n_{\textnormal{s}}^2,n_{0^+}n_{\textnormal{s}^+}^2, 1, 1, n_{\textnormal{ns}}, n_{\textnormal{ns}^+})(\FF_q).
\]
The entries list the genus $g$ of the curve, the size $q$ of the finite field,
the 4-tuple $(n_0,n_{0^+}, n_{\textnormal{ns}}, n_{\textnormal{ns}^+})$ and the number $\#X(\FF_q)$ of $\FF_q$-points of $X(n_0,n_{0^+}, n_{\textnormal{ns}}, n_{\textnormal{ns}^+})$.
When a lower bound $L_g(\FF_q)$ is not yet available, we list a curve of genus $g$ with $M$ points only if the corresponding upper bound satisfies $M_g(\FF_q)<q+1+\sqrt{2} (M-q-1)$, as it is done in the database \cite{manypoints} at the time of writing (September 2022).

\vskip 10 mm

\begin{center}
\begin{tabular}{c}
$
\begin{array}{c|c|c|c||c|c|c|c} 
g & q & (n_0,n_{0^+}, n_{\textnormal{ns}}, n_{\textnormal{ns}^+}) & \#X(\FF_q) & g & q & (n_0,n_{0^+}, n_{\textnormal{ns}}, n_{\textnormal{ns}^+}) & \#X(\FF_q) \\
%&\mathbb{F}_{11^5}  &\mathbb{F}_{13^5} &\mathbb{F}_{17^5} \\
\toprule
 5  &  5^2  & (1,572,1,1)  & 71 & 14 &   13^2  &  (7,19,1,3)  & 442\\%30
 6  &  3^2  & (1,398,1,1)  & 37 & 15 &   7^2  &  (1,956,1,1)  & 214\\

 7  &  11^5  & (1,12,1,7)  & 166589 & 15 &   13^5  &  (27,14,1,1)  & 384496\\
 8  &  13^2  & (1,4,9,1)  & 364  & 15 &   17^5  &  (1,80,3,1)  & 1445778\\
 9  &  3^5  & (5,17,2,1)  & 464 & 16 &   11^2  & (9,1,1,7)  & 414 \\
 9  &  7^5  & (99,1,1,1)  & 18968 & 16 &   11^5  & (3,4,1,7)  & 172432 \\
 9  &  11^3  & (4,43,1,1)  & 1812  & 17 &   13^2  &  (7,3,1,8)  & 500\\
 9  &  11^5  & (8,39,1,1)  & 167544  & 18 &   3^2  &  (1,878,1,1)  & 73\\
 9  &  13^5  & (96,1,1,1)  & 382096 & 18 &   7^2  &  (1,179,1,4)  & 224\\
 9  &  17^5  & (41,4,1,1)  & 1438108 & 18 &   17^2  &  (1,271,1,3)  & 746\\
 10  &  5^5  & (1,668,1,1)  & 4092 & 19 &   2^2  &  (225,1,1,1)  & 38\\%40
 10  &  11^5  & (3,104,1,1)  & 168744 & 19 &   13^2  &  (4,1,9,1)  & 582\\

 10  &  13^5  & (27,7,1,1)  & 382887 & 20 &   3^2  &  (1,1244,1,1)  & 80\\

 10  &  19^5  & (92,1,1,1)  & 2499156 & 20 &   17^2  &  (8,13,1,3)  & 862\\

 11  &   7^2  & (1,764,1,1)  & 176  & 21 &   11^2  &  (256,1,1,1)  & 464\\

 11  &   17^5  & (104,1,1,1)  & 1438748  & 21 &   17^2  &  (45,1,1,4)  & 824\\

 11  &   19^5  & (17,4,1,3)  & 2501908  & 22 &   2^2  &  (1,761,1,1)  & 43\\

 12  &   3^2  & (4,7,1,5)  & 58  & 22 &   3^2  &  (121,4,1,1)  & 80\\

 12 &   5^2  & (16,13,1,1)  & 122 & 22 &   13^2  &  (12,25,1,1)  & 594\\

 12 &   11^2  & (12,13,1,1)  & 338 & 23 &   5^2  & (1,887,1,1)  & 180\\

 12 &   7^2  & (1,718,1,1)  & 171 & 23 &   13^2  & (3,476,1,1)  & 594\\%50

 12 &   11^5  & (12,13,1,1)  & 170676 & 23 &   19^2  & (3,1,14,1)  & 988\\

 12 &   19^5  & (16,13,1,1)  & 1445778 & 24 &   5^2  & (1,412,1,3)  & 188\\

 13 &   5^2  & (4,143,1,1)  & 126 & 24 &   11^2  & (9,2,1,7)  & 498\\

 13 &   7^2  &  (1,599,1,1)  & 184 & 24 &   17^2  & (1,981,1,1)  & 880\\

 13 &   13^2  &  (9,1,1,8)  & 456 & 25 &   5^2  & (4,167,1,1)  & 204\\

 13 &   17^2  &  (144,1,1,1)  & 696 & 25 &   13^2  & (180,1,1,1)  & 672\\

 14 &   3^2  &  (1,734,1,1)  & 59 & 25 &   17^2  & (49,1,2,3)  & 990\\

 14 &   7^2  &  (1,734,1,1)  & 194 & 25 &   11^5  & (24,13,1,1)  & 180048\\%58

\end{array}
$
	\end{tabular}
	\captionof{table}{Improved bounds of the form $\#X(n_0,n_{0^+},n_{\mathrm{ns}},n_{\mathrm{ns}^+})(\FF_{q})$ for $n_0n_{0^+} n_{\textnormal{ns}}^2 n_{\textnormal{ns}^+}^2\le 10000$ and $g\leq 25$.}
 \label{tab:record1}
\end{center}

\newpage

\begin{center}
\begin{table}
\[
\begin{array}{c|c|c|c||c|c|c|c}
g & q & (n_0,n_{0^+}, n_{\textnormal{ns}}, n_{\textnormal{ns}^+}) & \#X(\FF_q) & g & q & (n_0,n_{0^+}, n_{\textnormal{ns}}, n_{\textnormal{ns}^+}) & \#X(\FF_q) \\
%&\mathbb{F}_{11^5}  &\mathbb{F}_{13^5} &\mathbb{F}_{17^5} \\
\toprule

 27 &   5^2  & (1,1509,1,1)  & 191 & 36 &   5^2  & (1,2327,1,1)  & 236\\

 27 &   11^2  & (9,76,1,1)  & 584 & 36 &   5^2  & (4,79,1,3)  & 243\\
 28 &   5^2  & (1,1336,1,1)  & 200 & 37 &   19^2  & (9,1,2,7) & 1452\\%30

 29 &   2^2  &  (1,1091,1,1)  & 55 & 38 &   3^2  & (1,1231,1,1)  & 131\\
 29 &   5^2  & (1,2004,1,1)  & 200 & 38 &   17^2  & (1,416,1,3)  & 1224\\
 29 &   17^2  & (99,4,1,1)  & 1000 & 39 &   5^2  & (1,1774,1,1)  & 260\\
 29 &   19^2  & (99,4,1,1)  & 1216 & 40 &   3^2  & (1,1756,1,1)  & 142\\
 30 &   3^2  &  (1,1375,1,1)  & 99 & 40 &   5^2  & (1,1559,1,1)  & 264\\
 30 &   13^2  &  (8,61,1,1)  & 730 & 41 &   5^2  & (3,83,1,4)  & 256\\
 30 &   19^2  &  (8,61,1,1)  & 1202 & 41 &   11^2  & (128,1,1,3)  & 880\\
 31 &   5^2  & (1,1532,1,1)  & 228 & 42 &   3^2  & (1,1279,1,1)  & 132\\
 31 &   13^2  & (81,7,1,1) & 744 & 42 &   5^2  & (1,2012,1,1)  & 296\\
 31 &   17^2  & (4,455,1,1)  & 1038 & 43 &   5^2  & (3,623,1,1)  & 266\\%40
 31 &   19^2  & (9,1,7,1)  & 1260 & 43 &   7^2  & (9,43,2,1)  & 444\\
 31 &   13^5  & (27,28,1,1)  & 398892 & 44 &   3^2  & (1,1966,1,1)  & 136\\
 32 &   3^2  & (1,1039,1,1)  & 112 & 44 &   5^2  & (1,1487,1,1)  & 289\\
 32 &   5^2  & (1,542,1,3)  & 213 & 45 &   3^2  & (1,1399,1,1)  & 145\\
 32 &   19^2  & (11,140,1,1)  & 1236 & 45 &   5^2  & (1,1427,1,1)  & 275\\
 33 &   5^2  & (1,1319,1,1)  & 225 & 46 &   5^2  & (4,263,1,1)  & 306\\
 34 &   3^2  & (1,1678,1,1)  & 113 & 47 &   2^2  & (1,2681,1,1)  & 74\\ 

 34 &   5^2  & (1,1223,1,1)  & 241 & 47 &   5^2  & (3,383,1,1)  & 298\\

 34 &   13^2  & (9,100,1,1)  & 798 & 48 &   5^2  & (1,796,1,3)  & 302\\
 34 &   17^2  & (16,1,1,7)  & 1260 & 48 &   13^2  & (9,97,1,1)  & 1058\\%50
 34 &   11^5  & (12,1,1,7)  & 180450 & 49 &   5^2  & (4,311,1,1)  & 315\\
 35 &   3^2  & (1,1916,1,1)  & 122 & 50 &   5^2  & (1,2396,1,1)  & 306\\
 35 &   5^2  & (1,1916,1,1)  & 242 & 50 &   7^2  & (1,2396,1,1)  & 506\\
 36 &   3^2  & (4,199,1,1)  & 130 & 50 &   19^2  & (9,146,1,1)  & 1750\\%54
\end{array}
\]
\caption{Improved bounds of the form $\#X(n_0,n_{0^+},n_{\mathrm{ns}},n_{\mathrm{ns}^+})(\FF_{q})$ for
$25<g\leq 50$.}\label{tab:record2}
\end{table}
\end{center}

\bibliographystyle{amsalpha}
\bibliography{ManyPoints.bbl}{}
\newpage

\color{black}
\section*{Appendix: More data}\label{sec:appendix}

In Tables \ref{tab1A}, \ref{tab2A}, \ref{tab3A}, \ref{tab4A}, \ref{tab5A}, \ref{tab6A}, \ref{tab7A} and \ref{tab8A}, we list the maximal values that we found for the number of points of the curves $X(n_0,n_{0^+}, n_{\textnormal{ns}}, n_{\textnormal{ns}^+})$ over $\mathbb{F}_{q}$, for $n_0n_{0^+} n_{\textnormal{ns}}^2 n_{\textnormal{ns}^+}^2\le 10000$ and $q=p^k$ a prime power with $p<20$, $k\le 5$ for $p$ odd and $k\leq 6$ for $p=2$. 
As in \Cref{sec:app} we restrict our search to the curves
$X(n_0,n_{0^+}, n_{\textnormal{ns}}, n_{\textnormal{ns}^+})$. The entries are of the type $(n_0,n_{0^+}, n_{\textnormal{ns}}, n_{\textnormal{ns}^+})\rightarrow \#X(n_0,n_{0^+}, n_{\textnormal{ns}}, n_{\textnormal{ns}^+})(\FF_q)$
and the field is indicated only once at the top of the column. All the values that improve previous known lower bounds $L_g(\FF_q)$ are in bold. When a lower bound $L_g(\FF_q)$ is not yet available, we display in bold curves with $M$ points only if the corresponding upper bound satisfies $M_g(\FF_q)<q+1+\sqrt{2} (M-q-1)$, as it is done in the database \cite{manypoints} at the time of writing (September 2022).
\newgeometry{text={190mm,230mm},centering}
\newpage
\tiny
\begin{center}

	\captionof{table}{Table for $\max \{\#X(n_0,n_{0^+},n_{\mathrm{ns}},n_{\mathrm{ns}^+})(\FF_q)\}$  with $q=2^k$ and $n_0n_{0^+} n_{\textnormal{ns}}^2 n_{\textnormal{ns}^+}^2\le 10000$ } \label{tab1A}
\end{center}

\newpage

\begin{table}
\[
\begin{array}{c|c|c|c|c|c}
g & \mathbb{F}_{3} & \mathbb{F}_{3^2} & \mathbb{F}_{3^3} &\mathbb{F}_{3^4}  &\mathbb{F}_{3^5} \\
\toprule
% q=3\\ Livello massimo 10000
  1  & (1,37,1,1)\to7  & (1,5,1,8)\to16  & (1,1,2,5)\to38  & (1,7,1,4)\to96  & (1,1,1,11)\to275   \\
 2  & (1,85,1,1)\to8  & (1,14,1,5)\to20  & (1,11,4,1)\to44  & (23,1,1,1)\to116  & (1,14,1,5)\to306   \\
 3  & (4,13,1,1)\to10  & (4,31,1,1)\to28  & (1,29,2,1)\to46  & (1,95,2,1)\to124  & (1,136,1,1)\to304   \\
 4  & (1,148,1,1)\to11  & (44,1,1,1)\to30  & (1,41,2,1)\to44  & (1,160,1,1)\to138  & (5,17,1,1)\to346   \\
 5  & (1,212,1,1)\to12  & (1,28,1,5)\to32  & (5,44,1,1)\to52  & (23,4,1,1)\to132  & (1,28,1,5)\to366   \\
 6  & (1,340,1,1)\to13  & \mathbf{(1,398,1,1)\to37}  & (1,8,5,1)\to56  & (1,272,1,1)\to142  & (121,1,1,1)\to364   \\
 7  & (1,296,1,1)\to13  & (4,71,1,1)\to39  & (31,4,1,1)\to54  & (61,2,1,1)\to137  & (4,95,1,1)\to366   \\
\end{array}
\]
\caption{Table for $\max\{\#X(n_0,n_{0^+},n_{\mathrm{ns}},n_{\mathrm{ns}^+})(\FF_q)\}$  with $q=3^k$ and
$n_0n_{0^+} n_{\textnormal{ns}}^2 n_{\textnormal{ns}^+}^2\le 10000$}
\label{tab2A}
\end{table}

\newpage

\begin{table}
\[
\begin{array}{c|c|c|c|c|c}
g & \mathbb{F}_{5} & \mathbb{F}_{5^2} & \mathbb{F}_{5^3} &\mathbb{F}_{5^4}  &\mathbb{F}_{5^5} \\
\toprule
% q=5\\ Livello massimo 10000
 1  & (1,7,1,4)\to10  & (1,1,1,21)\to36  & (1,132,1,1)\to148  & (1,1,1,11)\to675  & (1,13,1,4)\to3227   \\
 2  & (1,67,1,1)\to12  & (1,51,1,4)\to46  & (1,14,3,1)\to170  & (1,67,1,1)\to724  & (1,167,1,1)\to3328   \\
 3  & (4,19,1,1)\to15  & (1,19,1,12)\to56  & (1,28,3,1)\to192  & (1,29,2,1)\to738  & (13,1,1,4)\to3404   \\
 4  & (1,172,1,1)\to17  & (4,47,1,1)\to66  & (4,33,1,1)\to174  & (61,1,1,1)\to758  & (1,334,1,1)\to3530   \\
 5  & (1,278,1,1)\to16  & \mathbf{(1,572,1,1)\to71}  & (1,1,3,8)\to170  & (19,1,1,3)\to817  & (1,572,1,1)\to3631   \\
 6  & (1,163,1,1)\to19  & (8,13,1,1)\to74  & (3,53,1,1)\to178  & (9,13,1,1)\to902  & (12,7,1,1)\to3612   \\
 7  & (1,268,1,1)\to21  & (16,9,1,1)\to84  & (4,7,3,1)\to240  & (1,232,1,1)\to854  & (7,33,1,1)\to3644   \\
\end{array}
\]
\caption{Table for $\max\{\#X(n_0,n_{0^+},n_{\mathrm{ns}},n_{\mathrm{ns}^+})(\FF_q)\}$  with $q=5^k$ and
$n_0n_{0^+} n_{\textnormal{ns}}^2 n_{\textnormal{ns}^+}^2\le 10000$} \label{tab3A}
\end{table}

\newpage
\begin{table}
\[
\begin{array}{c|c|c|c|c|c}
g & \mathbb{F}_{7} & \mathbb{F}_{7^2} & \mathbb{F}_{7^3} &\mathbb{F}_{7^4}  &\mathbb{F}_{7^5} \\
\toprule
% q=7\\ Livello massimo 10000

 1  & (1,1,1,15)\to13  & (1,1,1,11)\to64  & (1,76,1,1)\to380  & (1,1,1,24)\to2496  & (1,4,1,5)\to17050   \\
 2  & (8,1,1,3)\to16  & (1,1,1,16)\to78  & (1,33,2,1)\to412  & (1,1,4,3)\to2590  & (1,8,1,5)\to17292   \\
 3  & (1,113,1,1)\to18  & (1,124,1,3)\to92  & (1,71,2,1)\to446  & (1,29,2,1)\to2684  & (1,290,1,1)\to17534   \\
 4  & (1,137,1,1)\to21  & (20,3,1,1)\to102  & (3,55,1,1)\to440  & (16,1,1,3)\to2778  & (1,22,3,1)\to17770   \\
 5  & (1,4,1,15)\to24  & (52,1,1,1)\to112  & (8,11,1,1)\to440  & (29,1,2,1)\to2808  & (11,9,1,1)\to18012   \\
 6  & (1,73,1,3)\to26  & (1,359,1,1)\to121  & (1,33,4,1)\to480  & (3,68,1,1)\to2774  & (4,3,1,5)\to17795   \\
 7  & (1,257,1,1)\to26  & (1,1,16,1)\to148  & (12,11,1,1)\to474  & (68,1,1,1)\to2844  & (8,1,1,5)\to18496   \\
\end{array}
\]
\caption{Table for $\max\{\#X(n_0,n_{0^+},n_{\mathrm{ns}},n_{\mathrm{ns}^+})(\FF_q)\}$  with $q=7^k$ and
$n_0n_{0^+} n_{\textnormal{ns}}^2 n_{\textnormal{ns}^+}^2\le 10000$}\label{tab4A}
\end{table}

\newpage
\begin{table}
\[
\begin{array}{c|c|c|c|c|c}
g & \mathbb{F}_{11} & \mathbb{F}_{11^2} & \mathbb{F}_{11^3} &\mathbb{F}_{11^4}  &\mathbb{F}_{11^5} \\
\toprule
% q=11\\ Livello massimo 10000
 1  & (1,20,1,3)\to18  & (1,1,1,15)\to144  & (1,43,1,1)\to1404  & (1,37,1,1)\to14875  & (1,5,1,12)\to161854   \\
 2  & (1,133,1,1)\to21  & (1,167,1,1)\to166  & (16,3,1,1)\to1468  & (1,161,1,1)\to15118  & (1,3,1,7)\to162656   \\
 3  & (1,109,1,1)\to25  & (1,105,1,4)\to188  & (4,19,1,1)\to1488  & (1,27,1,4)\to15287  & (1,312,1,1)\to163458   \\
 4  & (1,148,1,1)\to27  & (1,32,1,3)\to210  & (1,172,1,1)\to1600  & (3,1,1,8)\to15466  & (3,13,2,1)\to164260   \\
 5  & (1,181,1,1)\to32  & (1,1,2,15)\to232  & (19,3,2,1)\to1494  & (1,378,1,1)\to15699  & (3,1,1,7)\to165062   \\
 6  & (1,244,1,1)\to32  & (1,103,1,3)\to239  & (1,1,3,7)\to1622  & (12,7,1,1)\to15590  & (3,52,1,1)\to165864   \\
 7  & (1,229,1,1)\to37  & (1,64,1,3)\to276  & (4,5,3,1)\to1572  & (3,140,1,1)\to15636  & \mathbf{(1,12,1,7)\to166589}   \\
\end{array}
\]
\caption{Table for $\max\{\#X(n_0,n_{0^+},n_{\mathrm{ns}},n_{\mathrm{ns}^+})(\FF_q)\}$  with $q=11^k$ and
$n_0n_{0^+} n_{\textnormal{ns}}^2 n_{\textnormal{ns}^+}^2\le 10000$}\label{tab5A}
\end{table}

\newpage

\begin{table}
\[
\begin{array}{c|c|c|c|c|c}
g & \mathbb{F}_{13} & \mathbb{F}_{13^2} & \mathbb{F}_{13^3} &\mathbb{F}_{13^4}  &\mathbb{F}_{13^5} \\
\toprule
% q=13\\ Livello massimo 10000

1  & (1,1,1,21)\to21  & (1,1,1,11)\to196  & (1,4,1,5)\to2290  & (1,1,1,15)\to28899  & (1,1,1,24)\to372496   \\
 2  & (1,107,1,1)\to26  & (1,4,1,7)\to222  & (11,4,1,1)\to2382  & (1,165,1,1)\to29166  & (1,276,1,1)\to373698   \\
 3  & (1,43,1,3)\to29  & (1,4,1,9)\to242  & (1,49,1,3)\to2408  & (1,43,1,3)\to29573  & (1,37,1,3)\to374900   \\
 4  & (1,214,1,1)\to31  & (7,19,1,1)\to270  & (7,1,2,3)\to2478  & (1,172,1,1)\to29840  & (9,14,1,1)\to376094   \\
 5  & (1,323,1,1)\to34  & (1,23,2,3)\to291  & (5,44,1,1)\to2548  & (1,4,1,15)\to29950  & (32,3,1,1)\to377296   \\
 6  & (1,67,1,3)\to41  & (23,1,1,3)\to308  & (1,5,4,3)\to2478  & (1,609,1,1)\to29775  & (12,7,1,1)\to376094   \\
 7  & (9,7,2,1)\to36  & (9,16,1,1)\to324  & (1,5,8,1)\to2600  & (1,100,1,3)\to30426  & (9,16,1,1)\to377304   \\
\end{array}
\]
\caption{Table for $\max\{\#X(n_0,n_{0^+},n_{\mathrm{ns}},n_{\mathrm{ns}^+})(\FF_q)\}$  with $q=13^k$ and
$n_0n_{0^+} n_{\textnormal{ns}}^2 n_{\textnormal{ns}^+}^2\le 10000$} \label{tab6A}
\end{table}

\newpage
\begin{table}
\[
\begin{array}{c|c|c|c|c|c}
g & \mathbb{F}_{17} & \mathbb{F}_{17^2} & \mathbb{F}_{17^3} &\mathbb{F}_{17^4}  &\mathbb{F}_{17^5} \\
\toprule
% q=17\\ Livello massimo 10000
 1  & (1,190,1,1)\rightarrow25  & (1,1,1,11)\rightarrow324  & (1,61,1,1)\rightarrow5054  & (1,3,1,8)\rightarrow84096  & (1,2,1,9)\rightarrow1422141   \\
 2  & (1,23,1,3)\rightarrow30  & (1,161,1,1)\rightarrow358  & (3,29,1,1)\rightarrow5166  & (1,88,1,1)\rightarrow84670  & (1,14,1,5)\rightarrow1424424   \\
 3  & (1,127,1,1)\rightarrow36  & (1,71,2,1)\rightarrow392  & (4,1,1,5)\rightarrow5256  & (3,70,1,1)\rightarrow85244  & (33,1,1,1)\rightarrow1426584   \\
 4  & (16,5,1,1)\rightarrow42  & (1,31,2,3)\rightarrow426  & (7,20,1,1)\rightarrow5350  & (12,5,1,1)\rightarrow85818  & (1,7,5,1)\rightarrow1429090   \\
 5  & (1,316,1,1)\rightarrow38  & (9,1,4,1)\rightarrow448  & (5,31,1,1)\rightarrow5390  & (3,40,1,1)\rightarrow85496  & (1,208,1,1)\rightarrow1430904   \\
 6  & (1,92,1,3)\rightarrow46  & (3,2,5,1)\rightarrow460  & (1,47,1,4)\rightarrow5372  & (1,24,1,5)\rightarrow86345  & (8,13,1,1)\rightarrow1432818   \\
 7  & (1,3,1,16)\rightarrow40  & (16,9,1,1)\rightarrow516  & (5,11,1,3)\rightarrow5572  & (1,365,1,1)\rightarrow86577  & (1,129,2,1)\rightarrow1431052   \\
\end{array}
\]
\caption{Table for $\max\{\#X(n_0,n_{0^+},n_{\mathrm{ns}},n_{\mathrm{ns}^+})(\FF_q)\}$  with $q=17^k$ and
$n_0n_{0^+} n_{\textnormal{ns}}^2 n_{\textnormal{ns}^+}^2\le 10000$}\label{tab7A}
\end{table}

\newpage

\begin{table}
\[
\begin{array}{c|c|c|c|c|c}
g & \mathbb{F}_{19} & \mathbb{F}_{19^2} & \mathbb{F}_{19^3} &\mathbb{F}_{19^4}  &\mathbb{F}_{19^5} \\
\toprule
% q=19\\ Livello massimo 10000
 1  & (1,1,1,24)\rightarrow28  & (1,1,1,11)\rightarrow400  & (1,3,1,8)\rightarrow7024  & (1,11,1,3)\rightarrow131040  & (1,43,1,1)\rightarrow2478982   \\
 2  & (1,5,1,7)\rightarrow32  & (1,29,1,3)\rightarrow438  & (5,12,1,1)\rightarrow7188  & (1,5,1,7)\rightarrow131758  & (1,191,1,1)\rightarrow2482286   \\
 3  & (1,149,1,1)\rightarrow38  & (11,5,1,1)\rightarrow476  & (9,1,1,4)\rightarrow7352  & (1,23,1,4)\rightarrow132476  & (23,2,1,1)\rightarrow2484746   \\
 4  & (1,305,1,1)\rightarrow39  & (1,319,1,1)\rightarrow504  & (1,79,2,1)\rightarrow7249  & (1,260,1,1)\rightarrow133194  & (1,262,1,1)\rightarrow2487889   \\
 5  & (1,212,1,1)\rightarrow42  & (11,12,1,1)\rightarrow536  & (9,1,4,1)\rightarrow7352  & (1,130,1,3)\rightarrow132984  & (23,4,1,1)\rightarrow2490510   \\
 6  & (1,269,1,1)\rightarrow43  & (1,73,1,3)\rightarrow560  & (29,4,1,1)\rightarrow7420  & (8,13,1,1)\rightarrow133478  & (1,103,1,3)\rightarrow2490728   \\
 7  & (1,298,1,1)\rightarrow47  & (1,671,1,1)\rightarrow585  & (5,27,1,1)\rightarrow7464  & (3,35,2,1)\rightarrow133428  & (4,45,1,1)\rightarrow2490516   \\
\end{array}
\]
\caption{Table for $\max\{\#X(n_0,n_{0^+},n_{\mathrm{ns}},n_{\mathrm{ns}^+})(\FF_{q})\}$  with $q=19^k$ and
$n_0n_{0^+} n_{\textnormal{ns}}^2 n_{\textnormal{ns}^+}^2\le 10000$}\label{tab8A}
\end{table}

\end{document}